\documentclass[sn-mathphys,Numbered]{sn-jnl}

\usepackage{graphicx}%
\usepackage{multirow}%
\usepackage{amsmath,amssymb,amsfonts, mathtools}%
\usepackage{amsthm}%
\usepackage{mathrsfs}%
\usepackage[title]{appendix}%
\usepackage{xcolor}%
\usepackage{textcomp}%
\usepackage{manyfoot}%
\usepackage{booktabs}%
\usepackage{algorithm}%
\usepackage{algorithmicx}%
\usepackage{algpseudocode}%
\usepackage{listings}%

\theoremstyle{thmstyleone}%
\newtheorem{theorem}{Theorem}[section]
\newtheorem{proposition}[theorem]{Proposition}%
\newtheorem{lemma}[theorem]{Lemma}%
\newtheorem{corollary}[theorem]{Corollary}%

\newcommand{\X}{\mathcal{X}}
\newcommand{\Y}{\mathsf{Y}}

\newcommand{\com}{\mathsf{c}}

\newcommand{\cX}{\mathcal{X}}

\newcommand{\1}{\mathbf{1}}

\newcommand{\mykill}[1]{}

\newcommand{\SE}[1]{{\color{black}{#1}}}

\newcommand{\eins}{\1}%
\usepackage{geometry}
\begin{document}

\title[Hilbert's projective metric for functions of bounded growth and exponential convergence of Sinkhorn's algorithm]{Hilbert's projective metric for functions of bounded growth and exponential convergence of Sinkhorn's algorithm}


\author[1]{\fnm{Stephan} \sur{Eckstein}}\email{stephan.eckstein@uni-tuebingen.de}

\affil[1]{\orgdiv{Department of Mathematics}, \orgname{University of T\"{u}bingen}, \orgaddress{\street{Auf der Morgenstelle 10}, \city{T\"{u}bingen}, \postcode{72076}, \country{Germany}}}


\abstract{Motivated by the entropic optimal transport problem in unbounded settings, we study versions of Hilbert's projective metric for spaces of integrable functions of bounded growth. These versions of Hilbert's metric originate from cones which are relaxations of the cone of all non-negative functions, in the sense that they include all functions having non-negative integral values when multiplied with certain test functions. We show that kernel integral operators are contractions with respect to suitable specifications of such metrics even for kernels which are not bounded away from zero, provided that the decay to zero of the kernel is controlled. As an application to entropic optimal transport, we show exponential convergence of Sinkhorn's algorithm in settings where the marginal distributions have sufficiently light tails compared to the growth of the cost function.}




\maketitle
\section{Introduction}\label{sec:intro}
Hilbert's projective metric is a powerful geometric tool, particularly because many operators are contractions with respect to suitable specifications of this metric (see, e.g., \cite{birkhoff1957extensions,hilbert1895gerade,hyers1997topics,liverani1995decay,nussbaum1988hilbert}). 
A recent prominent example of such operators 
are the ones given by the Schr\"{o}dinger equations in entropic optimal transport, which famously lead to Sinkhorn's algorithm (see, e.g., \cite{ChenGeorgiouPavon.16b, Cuturi.13,DeligiannidisDeBortoliDoucet.21,Leonard.12,nutz2021introduction,SinkhornKnopp.67}). However, so far, the application of Hilbert's metric in this area is limited to settings in which the cost function of the entropic optimal transport problem is bounded on the relevant domain, thus excluding many cases of practical interest, like distance-based cost functions and unbounded marginal distributions. The goal of this paper is to remove such restrictions and allow for the applicability of Hilbert's metric in unbounded settings. Among others, this reveals that Sinkhorn's algorithm converges with exponential rate as soon as the marginal distributions have sufficiently light exponential tails compared to the order of growth of the cost function. Along the way, we establish that kernel integral operators are contractions with respect to suitable versions of Hilbert's metric, even if the kernel functions are not bounded away from zero.

\subsection{Setting and summary of the main results}\label{subsec:intromain}
Let $\X$ be a Polish space with metric $d_{\X}$, let $\mu, \nu \in \mathcal{P}(\X)$ be Borel probability measures on $\X$, and let $c : \X \times \X \rightarrow \mathbb{R}_+$ be a non-negative measurable function with \mbox{$\int c \,d\mu \otimes \nu < \infty$}, where $\mu\otimes\nu$ denotes the product measure of $\mu$ and $\nu$. The entropic optimal transport (EOT) problem is defined as (cf.~\cite{nutz2021introduction})
\begin{equation}\label{eq:eotintro}
\inf_{\pi \in \Pi(\mu, \nu)} \int c(x, y) \,\pi(dx, dy) + \varepsilon H(\pi, \mu \otimes \nu),\tag{EOT}
\end{equation}
where $\Pi(\mu, \nu) \subseteq \mathcal{P}(\X \times \X)$  is the set of all couplings with marginals $\mu$ and $\nu$, $\varepsilon$ is the penalization parameter and $H$ is the relative entropy. 

Sinkhorn's algorithm is a popular approach to numerically solve \eqref{eq:eotintro}. It is based on the fact that the unique optimizer $\pi^*$ of \eqref{eq:eotintro} is given by
\begin{equation}\label{eq:primalopti}
\frac{d\pi^*}{d \mu \otimes \nu}(x, y) = \exp\left(-\frac{c(x, y)}{\varepsilon}\right)\,g_1^*(x)\, g_2^*(y),
\end{equation}
where $g_1^*$ and $g_2^*$ solve the corresponding Schr\"{o}dinger equations
\begin{align*}
g_1^{*}(x) &= T_{\nu}(g_2^*)(x) := \left( \int \exp\left(-\frac{c(x, y)}{\varepsilon}\right) g_2^*(y) \,\nu(dy) \right)^{-1}, \\ 
g_2^{*}(y) &= T_{\mu}(g_1^*)(y) := \left( \int \exp\left(-\frac{c(x, y)}{\varepsilon}\right) g_1^*(x) \,\mu(dx) \right)^{-1}. 
\end{align*}
One way to state Sinkhorn's algorithm is to initialize $g_1^{(0)} = g_2^{(0)} \equiv 1$ and iteratively set, for $n \in \mathbb{N}$,
\begin{align}\label{eq:sinkhorniter}
\begin{split}
g_1^{(n+1)} &:= T_{\nu}(g_2^{(n)}),\\
g_2^{(n+1)} &:= T_{\mu}(g_1^{(n+1)}),\\
\frac{d\pi^{(n)}}{d\mu \otimes \nu}(x, y) &:= \exp\left(-\frac{c(x, y)}{\varepsilon}\right)\,g_1^{(n)}(x)\, g_2^{(n)}(y).
\end{split}
\end{align}
To understand this algorithm, the key building blocks of interest are kernel integral operators. In the following, we restrict our attention to the operator $T_{\mu}$, noting that the treatment of $T_{\nu}$ is analogous by symmetry. Consider the kernel integral operator $L_{K, \mu} : L^1(\mu) \rightarrow L^{\infty}(\mu)$ defined by
\[
L_{K, \mu}g(y) := \int K(x, y) g(x) \,\mu(dx) ~~~ \text{ for }g \in L^1(\mu),
\]
with $K: \X \times \X \rightarrow (0, 1]$ measurable. In relation to Sinkhorn's algorithm, we have $T_\mu(g_1) = (L_{K, \mu} g_1)^{-1}$ for $K(x, y) = \exp\Big(-\frac{c(x, y)}{\varepsilon}\Big)$.

To obtain exponential convergence of Sinkhorn's algorithm, one major part of the approach in this paper---in fact a major part in most prior works applying Hilbert's metric to the study of Sinkhorn's algorithm---is to show that kernel integral operators are contractions with respect to a suitable version of Hilbert's metric.
The benchmark in this regard is \cite[Theorem 2]{birkhoff1957extensions}, which establishes that, under the condition \mbox{$\inf_{x, y \in \X} K(x, y) > 0$}, there exists a version of Hilbert's metric with respect to which the operator $L_{K, \mu}$ is a contraction. However, for the specification of the kernel $K$ relevant for Sinkhorn's algorithm, this condition means that $\sup_{x, y \in \X} c(x, y) < \infty$, i.e., this limits the applicability to bounded cost functions. In what follows, we generalize the result in \cite{birkhoff1957extensions} and subsequently show the implications for Sinkhorn's algorithm.

We recall that each version of Hilbert's metric corresponds to a certain cone. For the study of Sinkhorn's algorithm in bounded settings, prior works use the cone of all non-negative functions, 
\[
C = \left\{g \in L^2(\mu) : g \geq 0 ~\mu\text{-a.s.}\right\}
\]
with corresponding Hilbert metric\footnote{Here, $\|\cdot\|_\infty$ denotes the $\mu$-essential supremum. For the more general definition of Hilbert's metric for arbitrary cones, we refer to Section \ref{subsec:hilbertmetric}.}
\[
d_{C}(g, \tilde{g}) = \log \|g/\tilde{g}\|_\infty + \log \|\tilde{g}/g\|_\infty, ~~~~g, \tilde{g} \in C.
\]
We observe that $d_{C}(g, \tilde{g})$ is infinite whenever either of the fractions $g/\tilde{g}$ or $\tilde{g}/g$ are unbounded. The key idea which we use to overcome this obstacle is based on the observation that we can rewrite
\begin{align*}
C &= \left\{g \in L^2(\mu) : \int f g \,d\mu \geq 0 \text{ for all } f \geq 0\right\},\\
d_{C}(g, \tilde{g}) &= \log \sup_{f \geq 0} \frac{\int fg \,d\mu}{\int f \tilde{g} \,d\mu} + \log \sup_{f \geq 0} \frac{\int f\tilde{g} \,d\mu}{\int f g \,d\mu}.
\end{align*}
This representation reveals that $C$ and $d_C$ are very particular cases of a more general construction, where the more general structure is based on using different classes of test functions $f \in \mathcal{F}$ instead of all non-negative functions $f \geq 0$. This means, we work with
\begin{align}\label{eq:conedef}
C_{\mathcal{F}} &= \left\{g \in L^2(\mu) : \int f g \,d\mu \geq 0 \text{ for all } f \in \mathcal{F} \right\},\\
d_{C_{\mathcal{F}}}(g, \tilde{g}) &= \log \sup_{f \in \mathcal{F}} \frac{\int fg \,d\mu}{\int f \tilde{g} \,d\mu} + \log \sup_{f \in \mathcal{F}} \frac{\int f\tilde{g} \,d\mu}{\int f g \,d\mu},\notag
\end{align}
for suitable classes of test functions $\mathcal{F}$. Roughly speaking, these test functions can only have a limited fraction of their mass in the tails, and hence $d_{C_{\mathcal{F}}}(g, \tilde{g})$ can be finite even if $g/\tilde{g}$ and $\tilde{g}/g$ are unbounded in the tails. \SE{We note that $C_{\mathcal{F}}$ is the dual cone to $\mathcal{F}$ and we refer to Section \ref{sec:newvariant} for the precise definition of the sets of test functions used.}

The first main result of this paper shows that with cones of the form $C_{\mathcal{F}}$, we can obtain contractivity of kernel integral operators even if the kernel is not bounded away from zero. To state the theorem, we fix $x_0 \in \X$ and set $B_r := \{x \in \X : d_{\X}(x, x_0) \leq r\}$ for $r > 0$.

\begin{theorem}\label{thm:introthmkernel}
	Assume that
	\begin{equation}\label{eq:thmintroass}
	\lim_{r \rightarrow \infty} \frac{\mu(\X \setminus B_r)}{\inf_{x, y \in B_r} K^2(x, y)} = 0.
	\end{equation}
	Then there exists a cone of the form $C_{\mathcal{F}}$ and $\kappa \in (0, 1)$ such that
	\[
	d_{C_{\mathcal{F}}}(L_{K, \mu}g, L_{K, \mu}\tilde{g}) \leq \kappa\,d_{C_{\mathcal{F}}}(g, \tilde{g}) ~~~\text{for all } g, \tilde{g} \in C_{\mathcal{F}}.
	\]
	\begin{proof}
		The result is a direct consequence of Corollary \ref{cor:cleankernel}.
	\end{proof}
\end{theorem}
Theorem \ref{thm:introthmkernel} generalizes the classical result from \cite{birkhoff1957extensions}, since the condition given in \eqref{eq:thmintroass} is easily seen to be satisfied if $\inf_{x, y \in \X} K(x, y) > 0$. \SE{As further examples, the theorem also applies in unbounded polynomial or exponential cases such as, for $p,\delta \in (0, \infty)$, 
	\begin{alignat*}{2}
		\inf_{x, y \in B_r} K^2(x, y) &\propto r^{-p} &\text{~~and~~} \mu(\X \setminus B_r) 
		&\propto r^{-p-\delta},\text{ or }\\
		\inf_{x, y \in B_r} K^2(x, y) &\propto \exp(-r^p) &\text{~~and~~} \mu(\X \setminus B_r) 
		&\propto \exp(-r^{p+\delta}),
	\end{alignat*}
where the exponential regime is the one which will be used for the analysis of Sinkhorn's algorithm.}

The second main result of the paper states that versions of Hilbert's metric corresponding to cones $C_{\mathcal{F}}$ can be used to obtain exponential convergence of Sinkhorn's algorithm in unbounded settings. Importantly, the obtained convergence not only holds with respect to this version of Hilbert's metric, but we also obtain exponential convergence with respect to more commonly used metrics like the total variation norm for the primal iterations. To state the theorem, let us write $C_{\mathcal{F}}^\mu = C_{\mathcal{F}}$ for cones as defined in \eqref{eq:conedef}, and $C_{\mathcal{F}}^\nu$ for cones of the same form, but for which $\mu$ is replaced by $\nu$.
\begin{theorem}\label{thm:introsinkhorn}
	Assume there exists $p > 0$ and $\delta > 0$ such that 
	\[
	\lim_{r \rightarrow \infty} \frac{\sup_{x, y \in B_r} c(x, y)}{r^p} = 0 ~~~\text{and}~~~\lim_{r \rightarrow \infty} \frac{\max\{\mu(\X \setminus B_r), \nu(\X \setminus B_r)\}}{\exp(-r^{p+\delta})} = 0.
	\]
	Then the iterations of Sinkhorn's algorithm converge with exponential rate to the optimizers of \eqref{eq:eotintro}. That is, there exist cones of the form $C_{\mathcal{F}_1}^\mu$ \hspace*{-0.8mm}and $C_{\mathcal{F}_2}^\nu$ and numbers $A > 0$ and $\kappa \in (0, 1)$ such that for all $n \in \mathbb{N}$,
	\begin{align*}
	d_{C_{\mathcal{F}_1}^\mu}\left(g_1^*, g_1^{(n)}\right) + d_{C_{\mathcal{F}_2}^\nu}\left(g_2^*, g_2^{(n)}\right) &\leq A \kappa^n, \\
	\big\|\pi^{(n)}-\pi^*\big\|_{TV} &\leq A \kappa^n.
	\end{align*}
	\begin{proof}
		The result follows from Theorem \ref{thm:sinkhorn} and Corollary \ref{cor:Sinkhorn}.
	\end{proof}
\end{theorem}
\mbox{Theorem \ref{thm:introsinkhorn}} is a generalization of results stated in bounded settings (cf.~\cite{ChenGeorgiouPavon.16b, Cuturi.13,DeligiannidisDeBortoliDoucet.21}) which require that $\sup_{x, y \in \X} c(x, y) < \infty$. A precise re-derivation for the case $\sup_{x, y \in \X} c(x, y) < \infty$ using Theorem \ref{thm:introsinkhorn} can be obtained by normalizing the metric $d_\X$ to $d_\X \leq 1$. \SE{A natural unbounded case covered by Theorem \ref{thm:introsinkhorn} is for instance given by $c(x, y) = \|x-y\|^q$ and $\mu$ and $\nu$ having Lebesgue densities which behave like $\exp(-\|x\|^{q+\delta})$ in the tails. In particular, the result applies to sub-Gaussian measures whenever the cost is of the form $\|x-y\|^{q}$ for $q < 2$, but not for $q \geq 2$. 
	
At this point, it is worth emphasizing that there are still large gaps towards fully understanding Sinkhorn's algorithm: While the assumptions of Theorem \ref{thm:introsinkhorn} (and \cite{conforti2023quantitative} discussed below) provide the most general known conditions for \emph{exponential convergence} of Sinkhorn's algorithm, there is still a gap both towards the weakest known condition for \emph{TV convergence without rate} ($\exp(\lambda c) \in L^1(\mu \otimes \nu)$ for some $\lambda > 0$, see \cite{NutzWiesel.22}), and the condition for \emph{existence of solutions} to \eqref{eq:eotintro} ($c \in L^1(\mu, \nu)$, see, e.g., \cite[Theorem 2.1]{nutz2021introduction}). Further, while Theorem \ref{thm:introsinkhorn} is stated in the most natural regime with polynomial growth of the cost function and exponential decay of the marginals, it is worth noting that for invertible functions $\varphi, \psi$, the optimal transport problem with marginals $\mu, \nu$ and cost $c$ is equivalent to the one with marginals $\varphi_{\#} \mu$, $\psi_{\#} \nu$ and cost $c \circ (\varphi^{-1}, \psi^{-1})$, and thus different regimes can be transformed to fit the setting as well.}

The remainder of this paper is structured as follows: Section \ref{subsec:related} gives an overview of the related literature. Section \ref{sec:notation} introduces the setting and notation, and recalls basic facts regarding Hilbert's metric. In Section \ref{sec:newvariant}, we introduce the particular versions of Hilbert's metric aimed at comparing integrable functions of bounded growth. Section \ref{sec:integraloperator} presents the results on kernel integral operators, and finally Section \ref{sec:Sinkhorn} gives the results related to Sinkhorn's algorithm.

\subsection{Related literature}\label{subsec:related}
While optimal transport theory has long been an important topic in probability theory and analysis (see, e.g., \cite{figalli2010mass,lott2009ricci,rachev1985monge,RachevRuschendorf.98a, RachevRuschendorf.98b,Villani.09}), it recently received increasing interest in the area of machine learning and related fields (see, e.g., \cite{WGAN.17, bunne2022proximal, caron2020unsupervised, CuturiPeyre.19,ge2021ota, SCHIEBINGER2019928, solomon2015convolutional}). One of the main drivers of this increased interest is the introduction of entropic regularization for the optimal transport problem in \cite{Cuturi.13}, which improves computational tractability by allowing for the use of Sinkhorn's algorithm \cite{Sinkhorn.64,SinkhornKnopp.67}. Sinkhorn's algorithm was initially called matrix scaling algorithm, and is also known as the iterative proportional fitting procedure \cite{DeligiannidisDeBortoliDoucet.21, Ruschendorf.95} and related to the DAD problem \cite{BorweinLewisNussbaum.94} and information projections \cite{Csiszar.75}.
Alongside this increase in applications, theoretical properties of entropic optimal transport have become increasingly well understood, for instance related to the approximation of optimal transport through its regularized version (see, e.g., \cite{AltschulerNilesWeedStromme.21,MR4505361,CarlierDuvalPeyreSchmitzer.17,carlier2023convergence,Chizat2020Faster,ConfortiTamanini.19,del2023improved,eckstein2023convergence, nenna2023convergence,nutz2022entropic}), statistical properties (see, e.g., \cite{GeneveyEtAl.16, mena2019statistical,pooladian2021entropic,rigollet2022sample}), and computational aspects (see, e.g., \cite{cuturi2014fast, CuturiPeyre.19,schmitzer2019stabilized,solomon2015convolutional}).

The convergence of Sinkhorn's algorithm in particular is studied in many recent works (see, e.g., \cite{Carlier.21, EcksteinNutz.21, ghosal2022convergence, NutzWiesel.22, Ruschendorf.95}), also using Hilbert's metric (see \cite{ChenGeorgiouPavon.16b, Cuturi.13,DeligiannidisDeBortoliDoucet.21}). All works using Hilbert's metric have, however, been restricted to bounded settings. 
While \cite{EcksteinNutz.21,ghosal2022convergence,NutzWiesel.22} establish convergence results for Sinkhorn's algorithm in quite high generality in unbounded settings, these results only yield qualitative convergence or polynomial rates. 
Only very recently, the string of literature \cite{conforti2022weak,conforti2023quantitative,greco2023non} introduces techniques which make it possible to obtain exponential rates of convergence for Sinkhorn's algorithm even in unbounded settings.
We emphasize that such results in unbounded settings are a major step forward in understanding Sinkhorn's algorithm, since numerical guarantees based on results which only apply in bounded settings naturally deteriorate for large scale problems. Results for unbounded settings can yield uniform rates and constants for arbitrarily large numerical problems. This is particularly relevant in cases in which entropic optimal transport problems are used as building blocks within an outer problem (see, e.g., \cite{pmlr-v119-ballu20a,cuturi2014fast,lei2022neural,yan2023learning,yang2023estimating}). Beyond numerical considerations, understanding Sinkhorn's algorithm in unbounded settings can yield insights on the entropic optimal transport functional in terms of stability properties (cf.~\cite{DeligiannidisDeBortoliDoucet.21}) and hence also for gradient flows (cf.~\cite{carlier2022lipschitz}).

The first and, so far, only result showing exponential convergence of Sinkhorn's algorithm in unbounded settings is given in \cite{conforti2023quantitative}. The methodology therein is mainly based on using techniques around differentiability, which naturally induces smoothness conditions on cost function and densities of marginal distributions.
The results we establish in this paper are thus based on completely different techniques compared to \cite{conforti2023quantitative} and lead to complementary results. 
To elaborate, the results in \cite{conforti2023quantitative} are stated for $\X = \mathbb{R}^d$ and quadratic cost function, $c(x, y) = \|x-y\|_2^2$, and lead to very strong convergence results including linear convergence of the derivatives of the iterations. While the techniques in \cite{conforti2023quantitative} appear applicable to more general twice continuously differentiable cost functions, the results derived in \mbox{Section \ref{sec:Sinkhorn}} do not require any smoothness assumptions on the cost function or marginals, but only controls on the tails. Further, the results in \cite{conforti2023quantitative} basically deal with a critical case where both $\inf_{x, y \in B_r} \exp(-c(x, y)) \propto \exp(-\lambda_1 r^2)$ and $\mu(\X \setminus B_r) \propto \exp(-\lambda_2 r^2)$ for some $\lambda_1, \lambda_2 > 0$. This can be seen as a critical case because both are of the same exponential rate, in which case the exponential convergence of Sinkhorn's algorithm can only be shown to hold for large enough penalization parameters. In this paper, we instead focus on non-critical cases where there is a bit of slack in terms of the exponential rate between the tails of the marginals and the decay of the kernel. This slack allows us to obtain an exponential rate of convergence for all penalization parameters, which is important since small penalization parameters are often the relevant regime in practice, particularly if the entropic problem is used as an approximation of the unregularized optimal transport problem. \SE{It is important to note that allowing for such slack excludes important cases, like the case of Gaussian (or sub-Gaussian) marginals and quadratic cost, where analytic solutions are known (see \cite{janati2020entropic}) and thus one might expect nice behavior of Sinkhorn's algorithm in such regimes as well. Further, while the results in the current paper apply for arbitrarily small regularization parameters, one must be alert to the fact that the contraction coefficient goes to one rapidly (roughly like $1-\exp(C/\varepsilon)$) for small regularization parameters, and thus the actual convergence result becomes correspondingly weaker in practice. While optimizing for constants would go beyond the scope of the current paper, we believe that a full understanding of Sinkhorn's algorithm must include an understanding of the constants, and we refer to \cite{chizat2024sharper} regarding an analysis of the optimal regime of the constants depending on the regularization parameter in bounded settings.}

Beyond the relation to optimal transport, Hilbert's metric (introduced in \cite{hilbert1895gerade}) is frequently applied to the study of asymptotic properties of operators (see, e.g.,  \cite{birkhoff1957extensions,hyers1997topics,kohlberg1982contraction,lemmens2012nonlinear,lemmens2013birkhoff,nussbaum1988hilbert}), often using a variety of different cones (see, e.g., \cite{bjorklund2010central, lemmens2019hilbert, liverani1995decay}). Particularly relevant in relation to the construction of cones we use in Section \ref{sec:newvariant} is the work in \cite{liverani1995decay}, where cones are constructed with the aim of obtaining contractivity of integral operators arising from dynamical systems. In contrast to the construction in Section \ref{sec:newvariant}, where we aim to obtain flexible tail behavior of functions in the cone, the construction in \cite{liverani1995decay} is rather aimed at smoothness properties of functions in the cone.

Lastly, kernel integral operators are fundamental objects in various areas, like kernel learning (see, e.g., \cite{steinwart2008support}) or integral equations and Fredholm theory (see, e.g., \cite{edmunds2018spectral, polyanin2008handbook}). Closely related to the topic of this paper is the study of asymptotic behavior of infinite-dimensional linear operators, for instance related to Perron-Frobenius theory (see, e.g., \cite{de1986irreducible,karlin1964existence,kohlberg1982contraction,lemmens2012nonlinear}), semigroups of positive operators (see, e.g., \cite{arendt1986one, marek1970frobenius}), or ergodic theory (see, e.g., \cite{bansaye2022non, nummelin2004general}). In relation to this literature, we emphasize that Theorem \ref{thm:introthmkernel} has obvious corollaries along the lines of a Perron-Frobenius theorem using Banach's fixed point theorem. 

\section{Setting and notation}\label{sec:notation}
Let $\cX$ be a Polish space with metric $d_\cX$. We always fix $x_0 \in \cX$ and define $B_r := \{x \in \cX : d_\cX(x_0, x) \leq r\}$ for $r > 0$. We further use the notation $B_r^{\com} := \cX \setminus B_r$ for the complement. 

Let $K: \cX \times \cX \rightarrow (0, 1]$ be measurable, and we call $K$ a kernel function (for kernel integral operators).
Functions $l : [0, \infty) \rightarrow [0, \infty)$ will always be assumed to be non-increasing, and we say that $K$ has lower decay $l$ if 
\[
\inf_{x, y \in B_r} K(x, y) \geq l(r)~~~\text{ for all } r > 0.
\]
For any function $f: \X \times \X \rightarrow \mathbb{R}$, its conjugate $f^\top: \X \times \X \rightarrow \mathbb{R}$ is defined by $f^\top(x, y) := f(y, x)$ for $x, y \in \X$.

Let $\mathcal{P}(\X)$ be the space of Borel probability measures on $\X$ and let $\mu \in \mathcal{P}(\X)$. Let $L^0$ be the set of measurable functions mapping $\X$ to $\mathbb{R}$, and $L^p(\mu) \subseteq L^0$ the subset of functions $f$ satisfying $\int |f|^p \,d\mu < \infty$ for $p \in [1, \infty)$, and $L^\infty(\mu) \subseteq L^0$  the subset of $\mu$-a.s.~bounded functions. Elements of $L^p(\mu)$ are as usual considered up to $\mu$-a.s.~equality. For $f \in L^0$, denote by $f_+ := \max\{0, f\}, f_- := -\min\{0, f\}$, so that $f = f_+ - f_-$.

Let $L_{K, \mu}: L^1(\mu) \rightarrow L^\infty(\mu)$ be defined by
\[
(L_{K, \mu} f) (x) := \int K(x, y) f(y) \,\mu(dy) \text{ for } f \in L^1(\mu) \text{ and } x \in \X.
\]

For $A \subseteq \X$, we always use the notation
\[
\eins_{A}(x) := \begin{cases}
1, & x \in A,\\
0, &\text{else,}
\end{cases}~~  \text{ for }x \in \X.
\]

\subsection{Hilbert's metric}\label{subsec:hilbertmetric}
This section shortly recalls the basics regarding Hilbert's metric. For more details, we refer for instance to \cite{birkhoff1957extensions, kohlberg1982contraction, lemmens2012nonlinear, lemmens2013birkhoff, nussbaum1988hilbert}.

For a Banach space $\Y$, we call $C \subseteq \Y$ a cone, if $C$ is closed, convex, $\lambda C := \{\lambda x : x \in C\} \subseteq C$ for all $\lambda > 0$ and $C \cap -C = \{0\}$. In later chapters, we will usually use $\Y = L^2(\mu)$. We note that the study of Hilbert's metric allows for more general settings---particularly regarding the Banach space structure on $\Y$ and closedness of $C$---but since the stated assumptions are imposed in parts of the literature and also satisfied in later chapters, we may as well impose them now.

We define the partial ordering $\leq_C$ by $f \leq_C g$ if $g-f \in C$ for $f, g \in \Y$.
For $f, g \in C$, we say that $g$ dominates $f$, if there exists $b > 0$ such that $f \leq_C b g$. We say $f \sim_C g$ if $f$ dominates $g$ and $g$ dominates $f$, in which case we define
\begin{align*}
M(f, g) &:= \inf\{ b > 0 : b g - f \in C\}, \\
m(f, g) &:= \sup\{ a > 0 : f - a g \in C\}.
\end{align*}
Since we assume closedness of $C$ the infimum and supremum are attained, if they are finite.
For $f, g \in C$, Hilbert's metric $d_C$ corresponding to $C$ is defined by 
\[
d_C(f, g) := \log\left(\frac{M(f, g)}{m(f, g)}\right) ~~~ \text{for } f \sim_C g, g\neq 0
\]
and $d_C(f, g) = \infty$, otherwise. \SE{This general definition will be precisely linked to the one of $d_{C_\mathcal{F}}$ given by \eqref{eq:conedef} in Lemma \ref{lem:calcdist}.}  We emphasize that $d_C$ indeed defines a projective metric on $C$, i.e., it satisfies all properties of a pseudo-metric and additionally $d_C(x, y) = 0$ implies $y = \lambda x$ for some $\lambda > 0$ (see, e.g., \cite[Lemma 2.1]{lemmens2013birkhoff}).

The following result, due to Birkhoff, is one of the key results which makes working with Hilbert's metric attractive.
\begin{theorem}[Birkhoff's contraction theorem \cite{birkhoff1957extensions}]\label{thm:birkhoff}
	Let $\Y_1, \Y_2$ be two Banach spaces and $C_1 \subseteq \Y_1, C_2 \subseteq \Y_2$ two cones.
	Let $T : \Y_1 \rightarrow \Y_2$ be a linear operator such that $T(C_1) \subseteq C_2$. Then $d_{C_2}(Tf, Tg) \leq d_{C_1}(f, g)$ for all $f, g \in C_1$. If further
	\[
	\Delta := \sup \{d_{C_2}(Tf, Tg) : f, g \in C_1, Tf\sim_{C_2} Tg \} < \infty, 
	\]
	then
	\[
	d_{C_2}(Tf, Tg) \leq \tanh(\Delta/4) ~d_{C_1}(f, g).
	\]
	\begin{proof}
		A modern proof is given, for instance, in \cite[Theorem A.4.1]{lemmens2012nonlinear}.
	\end{proof}
\end{theorem}

It is a natural question to consider which other distances can be bounded using Hilbert's metric for a certain cone. In this regard, the result below gives the main tool which will be used later on in Section \ref{sec:newvariant} to prove \mbox{Proposition \ref{prop:relationmetrics}}. It is a slightly adjusted version of \cite[Lemma 1.3]{liverani1995decay}.\footnote{\SE{Compared to \cite[Lemma 1.3]{liverani1995decay}, we relax the condition that $\|\cdot\|$ must be a norm and instead add the required properties on $\|\cdot\|$ as assumptions. Aside from the norm property, the assumptions of Lemma \ref{lem:boundothernorm} allow for slightly more flexibility regarding the constant $L$ which implicitly always equals one in \cite[Lemma 1.3]{liverani1995decay} whereas for instance in Proposition \ref{prop:relationmetrics} (ii), we need $L=3$.}}
\begin{lemma}\label{lem:boundothernorm}
	Let $C\subseteq \Y$ be a cone and $\|\cdot\| : \Y \rightarrow [0, \infty]$.\footnote{Here, $\|\cdot\|$ is not (necessarily) the norm on $\Y$, but it can be an arbitrary functional.} Let $f, g \in C$ with $\|f\| = \|g\| \neq 0$. Assume that
	\begin{align*}
	(g \geq_C a f) &\Rightarrow (\|g\| \geq a \|f\|) ~~~&&\text{for all }a>0,\\
	(b f \geq_C g) &\Rightarrow (b \|f\|\geq \|g\|) ~~~&&\text{for all }b>0,\\
	\exists L > 0: (\delta f \geq_C g-f \geq_C-\delta f) &\Rightarrow (\|g-f\| \leq L \delta \|f\|) ~~~&&\text{for all }\delta>0.
	\end{align*}
	Then
	\[
	\|g-f\| \leq L \|f\| \left(e^{d_{C}(f, g)} - 1\right).
	\]
	\begin{proof}
		Without loss of generality we assume $f \sim_C g$, otherwise the statement is trivially satisfied. Then, we can write $d_C(f, g) = \log(b/a)$ for some $b \geq a > 0$, where $b = M(f, g)$ and $a = m(f, g)$. From the first assumption and the definition of $a$, we find $\|g\|\geq a \|f\|$ and thus\SE{, using $\|g\|=\|f\|>0$, we get} $a \leq 1$ and from the second assumption and the definition of $b$, we obtain $\|g\|\leq b\|f\|$ and thus $b \geq 1$. We get
		\begin{align*}
		(b-a)f \geq_C (b-1) f \geq_C g-f \geq_C (a-1) f \geq_C -(b-a) f,
		\end{align*}
		which we can combine with the third assumption (using $\delta = (b-a)$ therein) to find $\|g-f\|\leq L(b-a) \|f\|$ and hence
		\[
		\|g-f\| \leq L \|f\| (b-a) \leq L \|f\|\frac{b-a}{a} = L \|f\| \left(e^{d_C(f, g)}-1\right),
		\]
		which completes the proof.
	\end{proof}
\end{lemma}

\section{Hilbert's metric for functions of bounded growth}\label{sec:newvariant}
This section introduces the new versions of Hilbert's metric to compare integrable functions of bounded growth. First, Subsection \ref{subsec:hilbertmotivation} gives an intuitive introduction to the family of cones which is used. In Subsection \ref{subsec:defandprop}, the formal definition for the cones is given and basic results are derived. More precisely, Lemma \ref{lem:basicpropFG} establishes that the defined sets are really cones and Lemma \ref{lem:calcdist} provides a simple formula to calculate the resulting Hilbert metric, Proposition \ref{prop:relationmetrics} shows how the defined versions of Hilbert's metric can bound other frequently used distances like the uniform norm on compacts or the $L^1$-norm, and finally Lemma \ref{lem:tailboundsfg} shows how to control the tails for the sets of functions which are used.

\subsection{Motivation and construction}\label{subsec:hilbertmotivation}
The aim of this subsection is to motivate the definition of the family of cones used in the remainder of this paper. In terms of mathematical content, readers may skip forward to Subsection \ref{subsec:defandprop}, as the purpose in the current subsection is purely to build intuition and explain the idea behind the given definitions. To this end, we restrict our attention to $\cX = \mathbb{R}_+$. Nevertheless, we will see later that the derived construction naturally transfers to the general case. 

We wish to define a family of cones whose corresponding versions of Hilbert's metric are well suited to compare functions of bounded growth.
The basic idea is to use cones of the form
\[
\mathcal{G} := \left\{g \in L^2(\mu) : \int f g \,d\mu \geq 0 \text{ for all } f \in \mathcal{F} \right\}
\]
for certain test functions $\mathcal{F}$,
with corresponding Hilbert metric
\begin{equation}\label{eq:dGdef}
d_{\mathcal{G}}(g, \tilde{g}) = \log \sup_{f \in \mathcal{F}} \frac{\int fg \,d\mu}{\int f \tilde{g} \,d\mu} + \log \sup_{f \in \mathcal{F}} \frac{\int f\tilde{g} \,d\mu}{\int f g \,d\mu}.
\end{equation}

Recall the baseline version of Hilbert's metric, where $\mathcal{F} = \mathcal{G}$ are the sets of all $\mu$-a.s.~non-negative functions, and $d_{\mathcal{G}}(g, \tilde{g}) = \log \|g/\tilde{g}\|_\infty + \log \|\tilde{g}/g\|_\infty$.
Intuitively, the main problem with this baseline version of Hilbert's metric is that in unbounded settings, the distance between $g$ and $\tilde{g}$ is only finite if $g$ and $\tilde{g}$ have \emph{precisely} the same level of growth. This means, even if $g(x) = 1+x$ and $\tilde{g}(x) = 1+x^{1+\varepsilon}$, then $d_{\mathcal{G}}(g, \tilde{g}) = \infty$. To alleviate this behavior, we clearly need a version of Hilbert's metric which produces smaller values in these cases. Considering formula \eqref{eq:dGdef}, it is thus natural to consider taking smaller sets of test functions $\mathcal{F}$.\footnote{In general, already Hilbert (cf.~\cite{hilbert1895gerade}) noted that increasing the size of the cone decreases the values of the corresponding version of Hilbert's metric. In our case, the test functions $\mathcal{F}$ and the cone $\mathcal{G}$ are inversely related in terms of size, i.e., increasing the size of $\mathcal{F}$ decreases the size of $\mathcal{G}$, and decreasing the size of $\mathcal{F}$ increases the size of $\mathcal{G}$.}

\begin{figure}
	\begin{minipage}{0.5\textwidth}
		\includegraphics[width=1.1\textwidth]{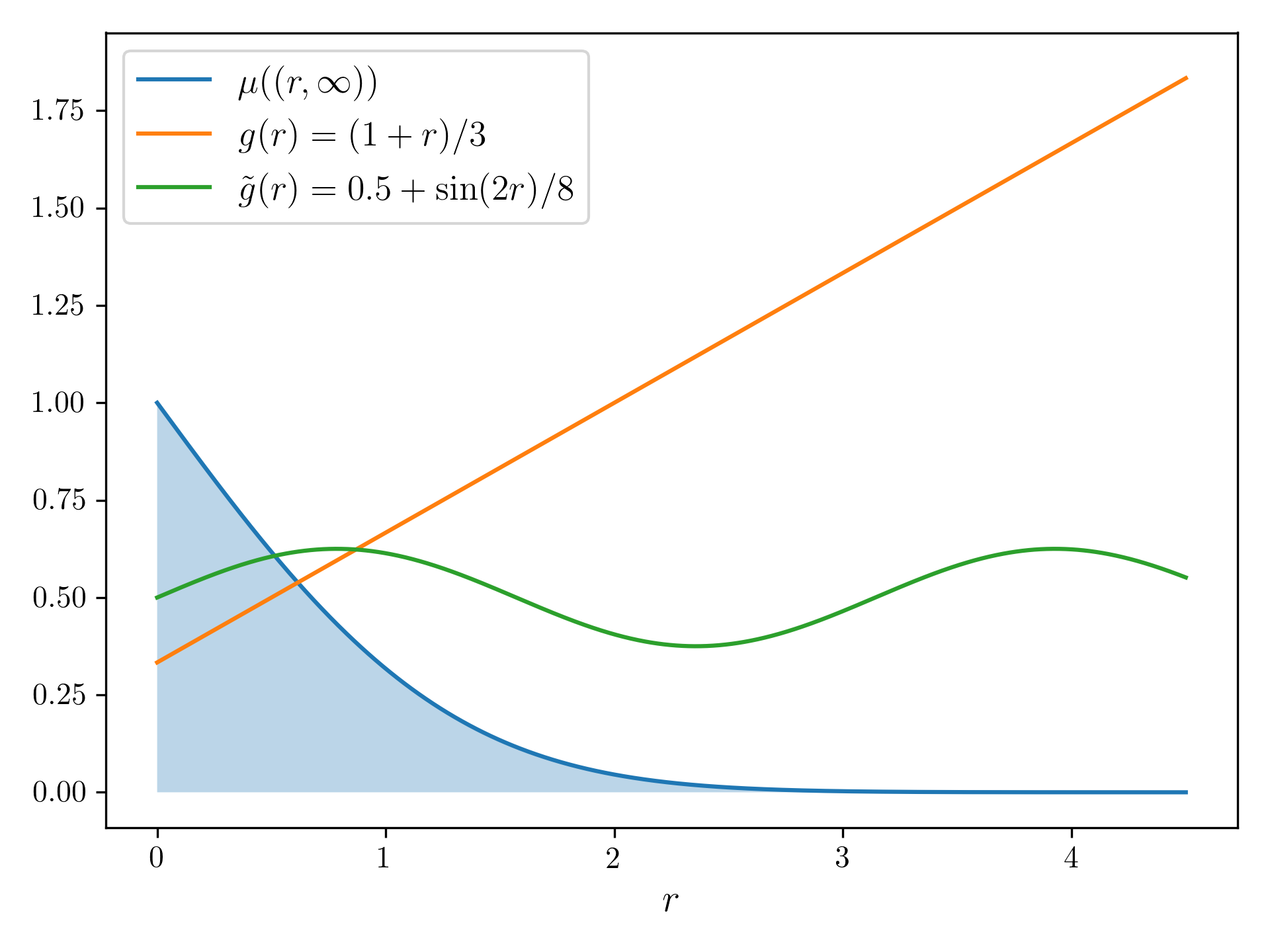}
	\end{minipage}%
	\begin{minipage}{0.5\textwidth}
		\includegraphics[width=1.1\textwidth]{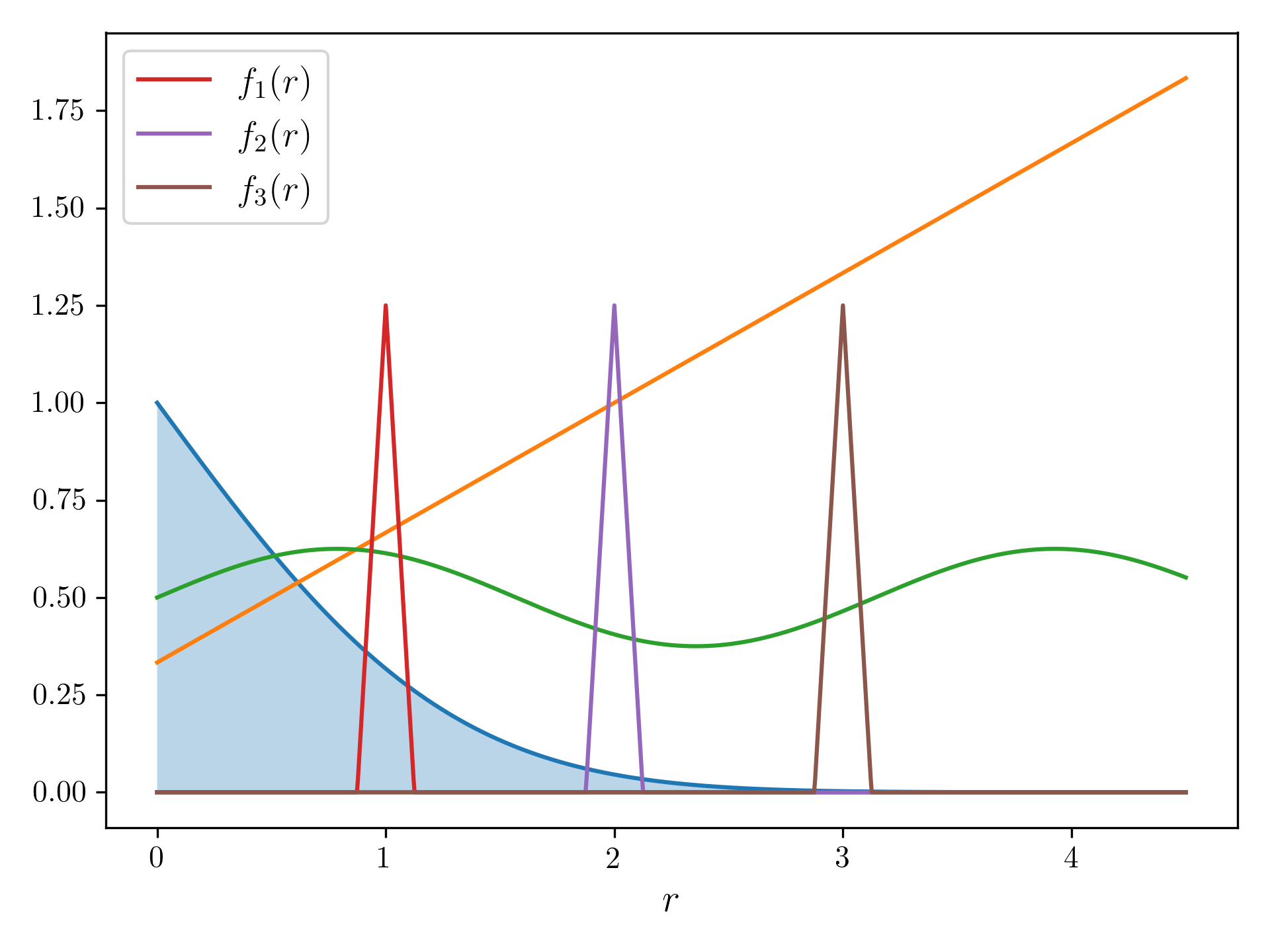}
	\end{minipage}
	\caption{The left hand side depicts a simple setting with a light-tailed measure, $\mu \sim |X|$, $X$ standard normal, and we wish to compare two functions of bounded growth, $g$ and $\tilde{g}$. The right hand image shows the start of a sequence of functions $f_1, f_2, f_3, \dots$ which show that formula \eqref{eq:dGdef} for the baseline version of Hilbert's metric evaluates to $d_{\mathcal{G}}(g, \tilde{g}) = \infty$.}\label{fig:example1}
\end{figure}

Figure \ref{fig:example1} showcases a simple case where we wish to compare two functions of bounded growth, but the baseline version of Hilbert's metric leads to a value of infinity. A sequence of test functions $f_1, f_2, f_3, \dots$ leading to this value in formula \eqref{eq:dGdef} is depicted on the right hand side of Figure \ref{fig:example1}. For such a sequence of functions, we say that $f_n$ has all of its mass in the ``tail'', i.e., on the interval $(n-1, \infty)$. This means, $\int f_n \,d\mu = \int_{(n-1, \infty)} f_n \,d\mu$. For a version of Hilbert's metric which produces a finite value in this example, we thus wish to impose conditions on the set of test functions which restricts such a behavior where mass is moved further and further into the tails. The basic idea to achieve this is to impose, for a fixed function \mbox{$\alpha: \mathbb{R}_+ \rightarrow \mathbb{R}_+$}, the condition
\begin{equation}\label{eq:conditionintuition}
\int_{(r, \infty)} f \,d\mu \leq \alpha(r) \int f \,d\mu ~ \text{ for all } r > 0.
\end{equation}
\begin{figure}
	\begin{minipage}{0.5\textwidth}
		\includegraphics[width=1.1\textwidth]{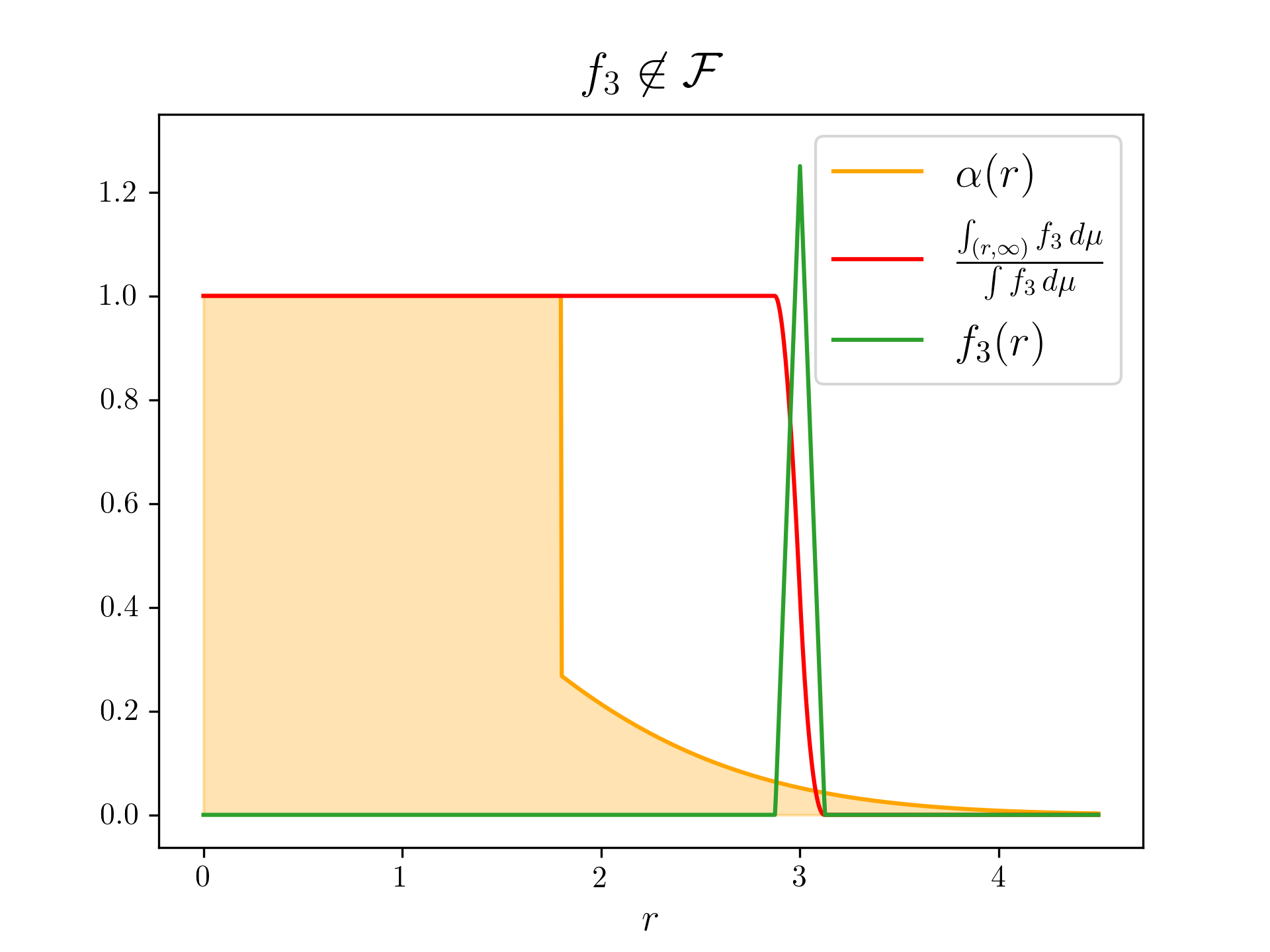}
	\end{minipage}%
	\begin{minipage}{0.5\textwidth}
		\includegraphics[width=1.1\textwidth]{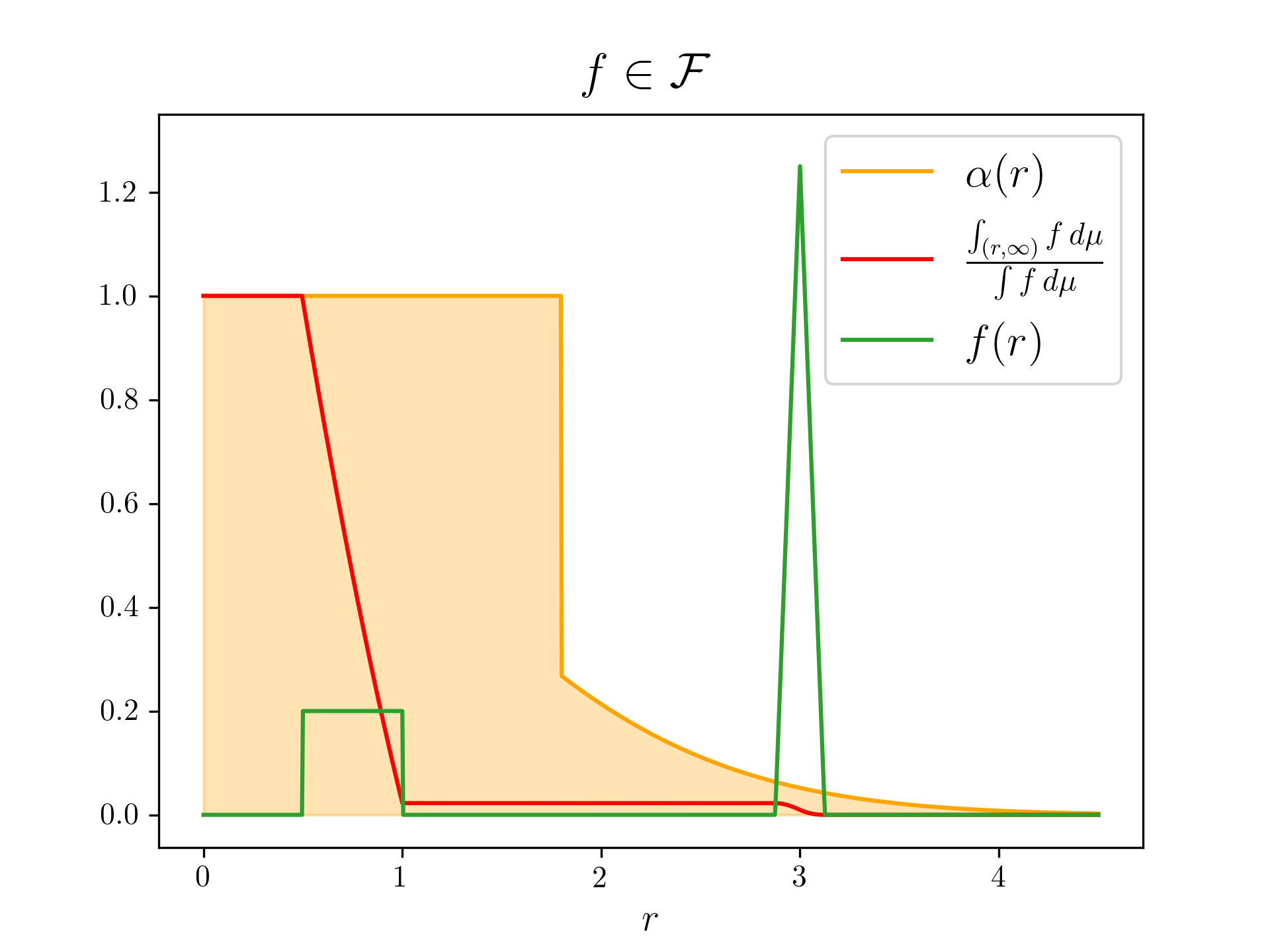}
	\end{minipage}
	\caption{The left hand side shows how condition \eqref{eq:conditionintuition} alleviates the undesired behavior in the example of Figure \ref{fig:example1}. With condition \eqref{eq:conditionintuition}, the sequence of functions shown in the right panel in Figure \ref{fig:example1} is not a valid sequence of test functions, since already $f_3$ puts a larger fraction of its mass into the tails than allowed by the function $\alpha$. The right hand side shows a function $f$ satisfying condition \eqref{eq:conditionintuition}. \SE{This function illustrates that condition \eqref{eq:conditionintuition} still allows test function to have positive values in the tails, as long as those are not the only positive values or disproportionately large.}}\label{fig:condition}
\end{figure}%
This means, for $f \in \mathcal{F}$ to hold, the function $f$ can at most put a fraction $\alpha(r)$ of its mass into the tail $(r, \infty)$. This condition is visualized in Figure \ref{fig:condition}. The left hand side of Figure \ref{fig:condition} shows how condition \eqref{eq:conditionintuition} excludes the sequence of test functions from Figure \ref{fig:example1}. The right hand side shows that test functions satisfying \eqref{eq:conditionintuition} can still \SE{have positive values in the tails, as long as they also have sufficient positive values in other parts of the space.}

In terms of the set $\mathcal{G}$, the given condition \eqref{eq:conditionintuition} allows for functions $g \in \mathcal{G}$ to take negative values in the tails, as long as those values are sufficiently small compared to the positive values of $g$ on other parts of the space. This means, while functions $g \in \mathcal{G}$ can take negative values, those negative values are nicely controlled. 

So far, our reasoning suggests that working with the set of test functions 
\[
\mathcal{F} = \left\{ f \in L^2(\mu) : f \geq 0~\mu\text{-a.s.~and } \int_{(r, \infty)} f \,d\mu \leq \alpha(r) \int f \,d\mu \text{ f.a.~}r > 0\right\}
\]
leads to a suitable version of the cone $\mathcal{G}$ and the corresponding version of Hilbert's metric $d_{\mathcal{G}}$. While this is mostly the case, one final problem which we need to alleviate is that with this construction, the positive values of functions $g \in \mathcal{G}$ are still completely unconstrained. This will particularly be important when we wish to calculate contraction coefficients for Theorem \ref{thm:birkhoff}, which includes taking a supremum over all functions in the cone $\mathcal{G}$. To control this supremum, we need a control also on the positive tails of functions in $\mathcal{G}$. 

This line of reasoning suggests that we allow for the set $\mathcal{F}$ to include functions which are partly negative. Indeed, so far, we can think of the positive parts of the test functions $f \in \mathcal{F}$ restricting the negative values of functions $g \in \mathcal{G}$. And in the same way, allowing for negative values of test functions $f \in \mathcal{F}$ can restrict the positive values of $g \in \mathcal{G}$. Of course, the cone $\mathcal{G}$ should still be close to the cone of all non-negative functions, so we only want positive values restricted in the tails. This line of reasoning motivates, for $m > 0$ and functions $\alpha, \tilde\alpha : \mathbb{R}_+ \rightarrow \mathbb{R}_+$, to use the set of test functions
\begin{align*}
\mathcal{F} := \Big\{f \in L^2(\mu) :~ &f \eins_{[0, m]} \geq 0~\mu\text{-a.s.}, \\
&\int_{(r, \infty)} f_+ \,d\mu \leq \alpha(r) \int f \,d\mu \text{ for } r \geq m, \\
&\int_{(r, \infty)} f_- \,d\mu \leq \tilde{\alpha}(r) \int f \,d\mu \text{ for } r \geq m
\Big\}
\end{align*}
This is precisely the definition of test functions we will adopt in the general setting in the following, where we merely replace intervals $[0, r]$ by balls $B_r$ and $(r, \infty)$ by $B_r^{\com}$ in arbitrary polish spaces.


%

\subsection{Definition and basic properties}
\label{subsec:defandprop}

\SE{We will always use a positive real number, $0 < m < \infty$}. We call a non-increasing function $\alpha : (0, \infty) \rightarrow [0, \infty)$ with $\lim_{r\rightarrow \infty} \alpha(r) = 0$ a tail function. Fix two tail functions $\alpha, \tilde{\alpha}$ and $\mu \in \mathcal{P}(\X)$. We make the standing assumptions that $\mu(B_{m}) > 0$ and $\alpha(r) > 0$ for all $r > 0$, which exclude degenerate cases. Define
\SE{\begin{align*}
\mathcal{F}^{\mu}_{\alpha, \tilde{\alpha}} := \Big\{f \in L^2(\mu) :~ &f \eins_{B_{m}} \geq 0, \\
&\int_{B_r^{\com}} f_+ \,d\mu \leq \alpha(r) \int f \,d\mu \text{ for } r \geq m, \\
&\int_{B_r^{\com}} f_- \,d\mu \leq \tilde{\alpha}(r) \int f \,d\mu \text{ for } r \geq m
\Big\}
\end{align*}}
as well as its dual cone
\begin{align*}
\mathcal{G}^{\mu}_{\alpha, \tilde{\alpha}} := \left\{g \in L^2(\mu) : \int fg \,d\mu \geq 0 \text{ f.a.~} f \in \mathcal{F}^{\mu}_{\alpha, \tilde{\alpha}}\right\}.
\end{align*}
Hereby, (in-)equalities for functions in $L^2(\mu)$ (like $f \eins_{B_{m}} \geq 0$) will always be understood in the $\mu$-a.s.~sense. \SE{For notational simplicity, we suppress the dependence of $\mathcal{F}^{\mu}_{\alpha, \tilde{\alpha}}$ and $\mathcal{G}^{\mu}_{\alpha, \tilde{\alpha}}$ on $m$, the reason being that $m$ will always be the same for all cones used within a given statement. Further, we note that the theory regarding the introduced versions of Hilbert's metric applies as well when using separate $m$ and $\tilde{m}$ for the inequalities for $f_+$ and $f_-$, respectively, in the definition of $\mathcal{F}^{\mu}_{\alpha, \tilde{\alpha}}$. This will, however, not be necessary for the main applications to Sinkhorn's algorithm and hence omitted for notational simplicity.}

While $\mathcal{F}^{\mu}_{\alpha, \tilde{\alpha}}$ is also a cone, the relevant cone which is used to define the versions of Hilbert's metric we are interested in is given by $\mathcal{G}^{\mu}_{\alpha, \tilde{\alpha}}$. That these sets are indeed cones is shown in the following.

\begin{lemma}\label{lem:basicpropFG}
	The sets $\mathcal{F}^{\mu}_{\alpha, \tilde{\alpha}}$ and $\mathcal{G}^{\mu}_{\alpha, \tilde{\alpha}}$ are cones.
	\begin{proof}		
		For $\mathcal{F}^{\mu}_{\alpha, \tilde{\alpha}}$, convexity, closedness and $\lambda \mathcal{F}^{\mu}_{\alpha, \tilde{\alpha}} \subseteq \mathcal{F}^{\mu}_{\alpha, \tilde{\alpha}}$ for $\lambda > 0$ are clearly true. Let further $f, -f \in \mathcal{F}^{\mu}_{\alpha, \tilde{\alpha}}$. We find $f\eins_{B_{m}} = 0$ and $\int f \,d\mu = 0$, and note that the latter implies $\int_{B_{m}^{\com}} |f| \,d\mu = 0$, which overall yields $f = 0$ (always understood in the $\mu$-a.s.~sense).
		
		For $\mathcal{G}^{\mu}_{\alpha, \tilde{\alpha}}$, again convexity and $\lambda \mathcal{G}^{\mu}_{\alpha, \tilde{\alpha}} \subseteq \mathcal{G}^{\mu}_{\alpha, \tilde{\alpha}}$ for $\lambda > 0$ are clearly true. Closedness is a simple consequence of the Cauchy-Schwartz inequality, since for $g_1, g_2, ... \in \mathcal{G}^{\mu}_{\alpha, \tilde{\alpha}}$ with $L^2(\mu)$-limit $g \in L^2(\mu)$ and $f \in \mathcal{F}^{\mu}_{\alpha, \tilde{\alpha}}$, we have
		\[
		\int fg \,d\mu \geq \liminf_{n\rightarrow \infty} \left[ \int f g_n \,d\mu - \|f\|_{L^2(\mu)} \|g-g_n\|_{L^2(\mu)} \right] \geq 0,
		\]
		which shows $g \in \mathcal{G}^{\mu}_{\alpha, \tilde{\alpha}}$.
		
		Finally, let $g, -g \in \mathcal{G}^{\mu}_{\alpha, \tilde{\alpha}}$. Then $\int gf \,d\mu = 0$ for all $f \in \mathcal{F}^{\mu}_{\alpha, \tilde{\alpha}}$, and thus $g \eins_{B_{m}} = 0$, since $f$ can be chosen as arbitrary non-negative functions on $B_{m}$. Further, with $f_{r, K} := K \eins_{B_{m}} + \eins_{\{g > 0\} \cap B_r}$, we find that $f_{r, K} \in \mathcal{F}_{\alpha, \tilde{\alpha}}^{\mu}$ for $K$ large enough since $\mu(B_{m}) > 0$ and $\alpha(r) > 0$. This implies $\int_{B_r} g_+ \,d\mu = 0$ and thus $\eins_{B_r} g_+ = 0$. We can similarly show that $\tilde{f}_{r, K} := K \eins_{B_{m}} + \eins_{\{g <0\} \cap B_r} \in \mathcal{F}^{\mu}_{\alpha, \tilde{\alpha}}$ for $K$ large enough, and thus $\eins_{B_r} g_- = 0$. Since both hold for all $r > 0$, we obtain $g=0$ and thus the claim.
	\end{proof}
\end{lemma}

The next result presents a formula for Hilbert's metric corresponding to the cone $\mathcal{G}^{\mu}_{\alpha, \tilde{\alpha}}$. This rigorously shows that the intuitive generalization of the baseline version of Hilbert's metric which we used in the introduction and \mbox{Subsection \ref{subsec:hilbertmotivation}} is indeed correct.
\begin{lemma}\label{lem:calcdist}
	For $g, \tilde{g} \in \mathcal{G}^{\mu}_{\alpha, \tilde{\alpha}}$, we have
	\[
	d_{\mathcal{G}^{\mu}_{\alpha, \tilde{\alpha}}}(g, \tilde{g}) = \log\left( \sup_{f \in \mathcal{F}^{\mu}_{\alpha, \tilde{\alpha}}} \frac{\int fg \,d\mu}{\int f\tilde{g} \,d\mu} \right) + \log\left( \sup_{f \in \mathcal{F}^{\mu}_{\alpha, \tilde{\alpha}}} \frac{\int f\tilde{g} \,d\mu}{\int fg \,d\mu} \right).
	\]
	Hereby, we use the convention $0/0 = 0$ and $a/0 = \infty$ for $a > 0$.
	\begin{proof}
		Let us first consider the case $g \sim_{\mathcal{G}^{\mu}_{\alpha, \tilde{\alpha}}} \tilde{g}$. Note that then
		\[
		\left(\int fg \,d\mu > 0 \Leftrightarrow \int f\tilde{g} \,d\mu > 0 \right) ~~~ \text{ for all } f \in \mathcal{F}_{\alpha, \tilde{\alpha}}^\mu.
		\]
		Thus, we can simply exclude the cases where the integrals are 0, which is done with the convention $0/0 = 0$. As in the proof of Lemma \ref{lem:basicpropFG}, $g, \tilde{g} \not\in \{0\}$ also excludes the case in which the integrals are 0 for all $f \in \mathcal{F}^{\mu}_{\alpha, \tilde{\alpha}}$.
		We can hence calculate 
		\begin{align*}
		M(g, \tilde{g}) &= \inf\left\{b \geq 0 : b \int \tilde{g} f \,d\mu \geq \int gf  \,d\mu \text{ f.a.~}f \in\mathcal{F}_{\alpha, \tilde{\alpha}}^\mu\right\}\\
		&=\inf\left\{b \geq 0 : b \geq \sup_{f \in \mathcal{F}_{\alpha, \tilde{\alpha}}^\mu} \frac{\int g f \,d\mu}{\int \tilde{g} f \,d\mu}\right\}\\
		&=\sup_{f \in \mathcal{F}_{\alpha, \tilde{\alpha}}^\mu} \frac{\int g f \,d\mu}{\int \tilde{g} f \,d\mu}
		\end{align*}
		and analogously
		\begin{align*}
		m(g,\tilde{g}) = \inf_{f \in \mathcal{F}^{\mu}_{\alpha, \tilde{\alpha}}} \frac{\int g f \,d\mu}{\int \tilde{g} f \,d\mu}.
		\end{align*}
		This implies 
		\[
		-\log(m(g,\tilde{g})) = \log\left(\sup_{f \in \mathcal{F}^{\mu}_{\alpha, \tilde{\alpha}}} \frac{\int \tilde{g} f \,d\mu}{\int g f \,d\mu}\right),
		\] 
		and hence we obtain the statement in the case $g \sim_{\mathcal{G}^{\mu}_{\alpha, \tilde{\alpha}}} \tilde{g}$. 
		
		If $g \not\sim_{\mathcal{G}^{\mu}_{\alpha, \tilde{\alpha}}} \tilde{g}$, we have $d_{\mathcal{G}_{\alpha, \tilde{\alpha}}^{\mu}}(g, \tilde{g}) = \infty$ by definition. Let us say without loss of generality that for any $b > 0$, $g \not\leq_{\mathcal{G}_{\alpha, \tilde{\alpha}}^{\mu}} b \tilde{g}$, which implies that for all $b > 0$, there exists $f \in \mathcal{F}_{\alpha, \tilde{\alpha}}^\mu$ such that $b \int \tilde{g}f \,d\mu < \int g f \,d\mu$. Potentially using the convention $a/0 = \infty$, this implies
		\[
		\sup_{f \in \mathcal{F}^{\mu}_{\alpha, \tilde{\alpha}}} \frac{\int g f \,d\mu}{\int \tilde{g} f \,d\mu} = \infty,
		\]
		which shows the claim in the case that $g \not\sim_{\mathcal{G}^{\mu}_{\alpha, \tilde{\alpha}}} \tilde{g}$, completing the proof.
	\end{proof}
\end{lemma}

The result below shows that we can bound other metrics using the introduced version of Hilbert's metric, making concrete the generally available method from \mbox{Lemma \ref{lem:boundothernorm}}.

\newpage
\begin{proposition}\label{prop:relationmetrics}\phantom{A}
	\begin{itemize}
		\item[(i)] Assume $\mu(A) > 0$ for all non-empty open sets $A \subseteq B_m$. Define $\|\cdot\|_m$ by $\|f\|_m := \sup_{x \in B_m} |f|$ for continuous functions $f \in L^2(\mu)$. Then
		\begin{align*}
		\|g-\tilde{g}\|_m &\leq \left[\exp\left(d_{\mathcal{G}^{\mu}_{\alpha, \tilde{\alpha}}}(g, \tilde{g})\right) - 1\right], \text{ for all}\\ g, \tilde{g} &\in \mathcal{S}_1 := \left\{\hat{g} \in \mathcal{G}^{\mu}_{\alpha, \tilde{\alpha}} : \|\hat{g}\|_m=1, \hat{g} \text{ continuous}, \hat{g} \geq 0\right\}.
		\end{align*}
		\item[(ii)] Assume $\alpha(r) \geq 2 \mu(B_r^{\com})$ for all $r \geq m$. Then
		\begin{align*}
		\|g-\tilde{g}\|_{L^1(\mu)} &\leq 3 \left[\exp\left(d_{\mathcal{G}^{\mu}_{\alpha, \tilde{\alpha}}}(g, \tilde{g})\right) - 1\right], \text{ for all}\\ g, \tilde{g} &\in \mathcal{S}_2 := \left\{\hat{g} \in \mathcal{G}^{\mu}_{\alpha, \tilde{\alpha}} : \|\hat{g}\|_{L^1(\mu)}=1, \hat{g} \geq 0\right\}.
		\end{align*}
	\end{itemize}
	\begin{proof}
		For ease of notation, set $C := \mathcal{G}^{\mu}_{\alpha, \tilde{\alpha}}$ in this proof.
		
		(i):
		We want to apply Lemma \ref{lem:boundothernorm} with constant $L=1$ (for all functions $g, \tilde{g} \in \mathcal{S}_1$). 
		
		The first two assumptions of Lemma \ref{lem:boundothernorm} follow easily: Indeed, if $g \geq_C a \tilde{g}$, we have
		\[
		\int_{B_m} g f \,d\mu \geq \int_{B_m} a \tilde{g} f \,d\mu
		\]
		for all non-negative functions $f$, and hence by assumption on $\mu$ and continuity of $g$ and $\tilde{g}$, this implies $g(x) \geq a \tilde{g}(x)$ for all $x \in B_m$. Since $g(x) = |g(x)|$, $\tilde{g}(x) = |\tilde{g}(x)|$ for all $x \in B_m$, we obtain $\|g\|_m \geq a \|\tilde{g}\|_m$. The second assumption (that $b \tilde{g} \geq_C g$ implies $b \|\tilde{g}\|_m \geq \|g\|_m$) follows analogously. 
		
		We finally show that the third assumption of Lemma \ref{lem:boundothernorm} holds with $L = 1$. Indeed, analogously to the above, we first find that $\delta\tilde{g} \geq_C g - \tilde{g} \geq_C -\delta \tilde{g}$ implies $(1+\delta)\tilde{g}(x) \geq g(x) \geq (1-\delta)\tilde{g}(x)$ for all $x \in B_m$. This implies $(g-\tilde{g})(x) \leq \delta \tilde{g}(x)$ and $(\tilde{g}-g)(x) \leq \delta \tilde{g}(x)$, which yields $\|g-\tilde{g}\|_m \leq \delta \|\tilde{g}\|_m$ and hence the third assumption, completing the proof of (i).
		
		(ii): We again show applicability of Lemma \ref{lem:boundothernorm}, this time with constant $L=3$. For the first part of the proof, we use the constant test function $f \equiv 1$. This satisfies $f \in \mathcal{F}^{\mu}_{\alpha, \tilde{\alpha}}$ since $\mu(B_r^{\com}) \leq \alpha(r)$. The first two assumptions in Lemma \ref{lem:boundothernorm} thus hold noting $g, \tilde{g} \geq 0$, which implies $\|g\|_{L^1(\mu)} = \int g \,d\mu = \int fg \,d\mu$, and similarly for $\tilde{g}$. 
		
		We further show that the third assumption is satisfied with $L = 3$, i.e., for all $g, \tilde{g} \in \mathcal{S}_2$, we show that
		\[(\delta \tilde{g} \geq_C g-\tilde{g} \geq_C-\delta \tilde{g}) \Rightarrow (\|g-\tilde{g}\|_{L^1(\mu)} \leq 3 \delta \|\tilde{g}\|_{L^1(\mu)}) ~~~\text{f.a.~}\delta>0.\]
		To see that this is true, note that $(\delta \tilde{g} \geq_C g-\tilde{g} \geq_C-\delta \tilde{g})$ yields that for all $f_1, f_2 \in \mathcal{F}_{\alpha, \tilde{\alpha}}^{\mu}$,
		\begin{align*}
		\int f_1 ((1+\delta)\tilde{g} - g) \,d\mu &\geq 0,\\
		\int f_2 (g - (1-\delta) \tilde{g}) \,d\mu &\geq 0,
		\end{align*}
		and thus
		\begin{equation}\label{eq:f1f2}
		\int (\tilde{g}-g)(f_1-f_2) + \delta \tilde{g} (f_1 + f_2) \,d\mu \geq 0.
		\end{equation}
		We choose $f_1 = 1 + \eins_{g > \tilde{g}}$ and $f_2 = 1+\eins_{\tilde{g} > g}$. Since
		\[
		\int_{B_r^{\com}} f_i \,d\mu \leq 2 \mu(B_r^{\com}) \leq \alpha(r) \leq \alpha(r) \int f_i \,d\mu,
		\]
		we have $f_i \in \mathcal{F}^{\mu}_{\alpha, \tilde{\alpha}}$ and hence \eqref{eq:f1f2} yields
		\[
		0 \leq \int (\tilde{g}-g)(\eins_{g > \tilde{g}} - \eins_{\tilde{g} > g}) + 3 \delta \tilde{g} \,d\mu = - \|\tilde{g}-g\|_{L^1(\mu)} + 3 \delta \|\tilde{g}\|_{L^1(\mu)},
		\]
		which show that the third assumption of Lemma \ref{lem:boundothernorm} holds and thus completes the proof.
	\end{proof}
\end{proposition}

In the following, we will establish tail bounds for functions $g \in \mathcal{G}^{\mu}_{\alpha, \tilde{\alpha}}$. To this end, define
\begin{align}\label{eq:defbeta}
\begin{split}
{\tilde{\beta}}(r) &:= \frac{1}{\mu(B_{m})}\sup_{a \geq r}\frac{\mu(B_a^{\com})}{\alpha(a) \land 1} ~~ \text{ for } r \geq m, \\
\beta(r) &:= \frac{1}{\mu(B_{m})} (\sup_{a \geq r} \frac{\mu(B_a^{\com})}{\tilde{\alpha}(a)} + \mu(B_r^{\com})) ~~ \text{ for } r \geq m.
\end{split}
\end{align}
The functions $\beta$ and $\tilde{\beta}$ will control the tails of $g \in \mathcal{G}^{\mu}_{\alpha, \tilde{\alpha}}$ in a similar fashion as the functions $\alpha$ and $\tilde{\alpha}$ control the tails of $f \in \mathcal{F}^{\mu}_{\alpha, \tilde{\alpha}}$, as established Lemma \ref{lem:tailboundsfg} (i) below.\footnote{\SE{Note that a-priori, the functions $\beta$ and $\tilde{\beta}$ are allowed to attain the value $\infty$, which however clearly leads to empty statements in Lemma 3.4 below.}}  In terms of notation, we will be always be consistent in the sense that if different tail functions $\alpha_i, \tilde{\alpha}_i$, $i \in \mathbb{N}$, are used, then we also write $\beta_i, \tilde{\beta}_i$, $i \in \mathbb{N}$ for the corresponding functions defined in \eqref{eq:defbeta}. 

Part (ii) of the lemma below derives simple conclusions which are stated for the sets of functions normalized to have integral one, i.e., for
\[
\mathcal{F}_{\alpha, \tilde{\alpha}}^{\mu, norm} := \left\{f \in \mathcal{F}_{\alpha, \tilde{\alpha}}^\mu  : \int f \,d\mu = 1 \right\} ~ \text{ and } ~ 
\mathcal{G}_{\alpha, \tilde{\alpha}}^{\mu, norm} := \left\{g \in \mathcal{G}_{\alpha, \tilde{\alpha}}^\mu  : \int g \,d\mu = 1 \right\}.
\]
\begin{lemma}\label{lem:tailboundsfg}
	The following hold:
	\begin{itemize}
		\item[(i)]For $g \in \mathcal{G}_{\alpha, \tilde\alpha}^{\mu}$, we have
		\begin{align}\label{eq:tailneg}
		\int_{B_r^{\com}} g_- \,d\mu &\leq {\tilde{\beta}}(r) \int_{B_{m}} g \,d\mu ~\text{ for all } r \geq m,\\\int_{B_r^{\com}} g_+ \,d\mu &\leq \beta(r) \int_{B_{m}} g \,d\mu ~\text{ for all } r \geq m.\label{eq:tailpos}
		\end{align}
		\item[(ii)] For $f \in \mathcal{F}_{\alpha, \tilde{\alpha}}^{\mu, norm}$ and $g \in \mathcal{G}_{\alpha, \tilde{\alpha}}^{\mu, norm}$, we have
		\begin{align}
		\int_{B_r^{\com}} |f| \,d\mu &\leq \alpha(r) + \tilde{\alpha}(r), ~~~~~~~~~ \text{ for all } r \geq m,\label{eq:tailbound1}\\
		\int |f| \,d\mu &\leq 1+ 2 \tilde{\alpha}(m),\label{eq:tailbound2}\\
		\int_{B_r^{\com}} |g| \,d\mu &\leq \frac{\tilde{\beta}(r) + \beta(r)}{1-\tilde{\beta}(m)},~~~~~~~~\,\text{  for all } r \geq m,\label{eq:tailbound3}\\
		\int |g|\,d\mu &\leq 1+2\frac{\tilde{\beta}(m)}{1-\tilde{\beta}(m)}. \label{eq:tailbound4}
		\end{align}
	\end{itemize}
	
	\begin{proof}
		\begin{itemize}
			\item[(i)] First, we treat equation \eqref{eq:tailneg}. The inequality will follow by taking $f = u \eins_{B_{m}} + U \eins_{\{g<0\} \cap B_r^{\com}}$ for suitable $u, U \geq 0$. We normalize $\int f \,d\mu = 1$ and thus 
			\begin{align*}
			1 &= u \mu(B_{m}) + U \mu(\{g<0\} \cap B_r^{\com}) \\
			\Leftrightarrow ~u &= (1-U \mu(\{g<0\} \cap B_r^{\com}))/\mu(B_{m}).
			\end{align*}
			For $u \geq 0$ to hold, we need $U \mu(\{g<0\} \cap B_r^{\com}) \leq 1$.
			For $f \in \mathcal{F}_{\alpha, \tilde\alpha}^{\mu}$ to hold, noting that $f \geq 0$, we further only need $\int_{B_a^{\com}} f \,d\mu \leq \alpha(a)$ for all $a \geq m$. This reduces to $U \leq \inf_{a\geq r} \frac{\alpha(a)}{\mu(\{g<0\} \cap B_a^{\com})}$ (note that we can ignore $a < r$ since $f$ is zero on $B_{r}\setminus B_{m}$ and $a\mapsto \alpha(a)$ is non-increasing). Setting 
			\[U := \inf_{a\geq r} \frac{\alpha(a) \land 1}{\mu(B_a^{\com})} \leq \inf_{a\geq r} \frac{\alpha(a)}{\mu(\{g<0\} \cap B_a^{\com})}
			\]
			we thus obtain that $f \in \mathcal{F}_{\alpha, \tilde\alpha}^{\mu}$ and
			\[
			U \mu(\{g < 0\} \cap B_r^{\com}) \leq \frac{\alpha(r) \land 1}{\mu(B_r^{\com})} \mu(\{g<0\} \cap B_r^{\com}) \leq 1
			\]
			and hence $u \geq 0$. \SE{Notably, the case $U = 0$ corresponds to the case in which $\tilde{\beta} \equiv \infty$, in which case the statement of the current lemma is trivial. Thus, we can assume $U > 0$ in the following}		and further calculate
			\[
			\frac{u}{U} = \frac{1}{\mu(B_{m})} \left( \frac{1}{U} - \mu(\{g < 0\} \cap B_r^{\com}) \right) \leq \frac{1}{\mu(B_{m})} \frac{1}{U} = {\tilde{\beta}}(r)
			\]
			which yields equation \eqref{eq:tailneg} by noting that
			\[
			\int fg \,d\mu \geq 0 \Leftrightarrow u \int_{B_{m}} g \,d\mu \geq U \int_{B_r^{\com}} g_- \,d\mu.
			\]
			
			Next, we treat equation \eqref{eq:tailpos}. The proof is analogous to the one of equation \eqref{eq:tailneg} and follows by defining a function
			\[
			f:= u \eins_{B_{m}} - U \eins_{\{g>0\} \cap B_r^{\com}},
			\]
			with $u$ such that $\int f\,d\mu = 1$ and $U := \inf_{a\geq r} \frac{\tilde{\alpha}(a)}{\mu(B_a^{\com})}$, which yields \SE{$u = (1 + U \mu(\{g > 0\} \cap B_r^{\com}))/(\mu(B_m))$, $f \in \mathcal{F}_{\alpha, \tilde\alpha}^{\mu}$ and $u/U = \frac{1}{\mu(B_m)} (\sup_{a \geq r} \frac{\mu(B_a^{\com} \cap \{g > 0\})}{\tilde{\alpha}(a)} + \mu(B_r^{\com})) \leq \beta(r)$.}
			
			\item[(ii)] The bound in equation \eqref{eq:tailbound1} follows directly from the definition of $\mathcal{F}^{\mu}_{\alpha, \tilde{\alpha}}$ by writing $|f| = f_+ + f_-$. Bound \eqref{eq:tailbound2} follows by writing $|f| = f + 2 f_-$ and noting that $\int f_- \,d\mu = \int_{B_{m}^c} f_- \,d\mu$.
			
			For the bounds in equations \eqref{eq:tailbound3} and \eqref{eq:tailbound4}, we first find that \eqref{eq:tailneg} yields $\int_{B_{m}^{\com}} g_- \,d\mu \leq {\tilde{\beta}}(m) \int_{B_{m}} g \,d\mu$ and thus 
			\[\int g \,d\mu \geq \int_{B_{m}} g \,d\mu - \int_{B_{m}^{\com}} g_-\,d\mu \geq (1-{\tilde{\beta}}(m)) \int_{B_{m}} g \,d\mu.\] Hence, \eqref{eq:tailbound3} follows by writing $|g| = g_+ + g_-$ and using \eqref{eq:tailneg} and \eqref{eq:tailpos}. Finally, \eqref{eq:tailbound4} follows by writing $|g| = g + 2g_-$, using $\eqref{eq:tailneg}$ and noting $\int g_- \,d\mu = \int_{B_{m}^c} g_- \,d\mu$.\qedhere
		\end{itemize}
	\end{proof}
\end{lemma}

\section{Contractivity of kernel integral operators}\label{sec:integraloperator}

This section studies contractivity of kernel integral operators with respect to versions of Hilbert's metric as introduced in Section \ref{sec:newvariant}. First, \mbox{Proposition \ref{prop:nonexpansive}} establishes sufficient conditions so that kernel integral operators are, at least, non-expansive. The remainder of section then aims at establishing conditions so that kernel integral operators are strict contractions. The main result in this regard is given by Theorem \ref{thm:generalkernel}.

Let $0<m<\infty$, and $K : \X \times \X \rightarrow (0, 1]$ be a kernel function with lower decay $l : [0, \infty) \rightarrow [0, \infty)$, i.e., $\inf_{x, y \in B_r} K(x, y) \geq l(r)$ for all $r \geq 0$. Let $\mu, \nu \in \mathcal{P}(\X)$ and $\alpha_1, \tilde{\alpha}_1$ and $\alpha_2, \tilde{\alpha}_2$ be two pairs of tail functions. \SE{We recall that $m$ implicitly determines the region $B_m$ in which the sets of test functions are purely positive, and $m$ is always the same for all cones used.} The following equations give assumptions on the relationship between the decays of kernel, measures and tail functions, which are sufficient to establish non-expansiveness of kernel integral operators in Proposition \ref{prop:nonexpansive} below.
\begin{align}
l(m)(1-\alpha_1(m)) - \tilde{\alpha}_1(m) &\geq 0\label{eq:gengeqclosed}\\
\bar{a} := \sup_{r \geq m} (l(r) \nu(B_r) (1-\alpha_1(r)) - \tilde{\alpha}_1(m)) &> 0  \label{eq:genintegralcomparable}\\
\frac{1+\tilde{\alpha}_1(m)}{\bar{a}} \, \nu(B_r^{\com}) &\leq \alpha_2(r) ~~~(\text{for } r \geq m) \label{eq:gentailplusclosed}\\
\frac{\tilde{\alpha}_1(m)}{\bar{a}} \, \nu(B_r^{\com}) &\leq \tilde{\alpha}_2(r) ~~~(\text{for } r \geq m) \label{eq:gentailminusclosed}
\end{align}

\begin{proposition}\label{prop:nonexpansive}
	Assume \eqref{eq:gengeqclosed}-\eqref{eq:gentailminusclosed} hold. Then, we have
	\begin{align*}
	L_{K, \mu}(\mathcal{F}_{\alpha_1, \tilde{\alpha}_1}^{\mu}) \subseteq \mathcal{F}_{\alpha_2, \tilde{\alpha}_2}^{\nu} ~~ \text{ and thus } ~~L_{K^\top, \nu}(\mathcal{G}_{\alpha_2, \tilde{\alpha}_2}^{\nu}) \subseteq \mathcal{G}_{\alpha_1, \tilde{\alpha}_1}^{\mu}.
	\end{align*}
	\begin{proof}
		We shortly recall the normalization $\|K\|_\infty \leq 1$, which is used at many points throughout the proof.
		We first show $L_{K, \mu}(\mathcal{F}_{\alpha_1, \tilde{\alpha}_1}^{\mu}) \subseteq \mathcal{F}_{\alpha_2, \tilde{\alpha}_2}^{\nu}$. Let $f \in \mathcal{F}_{\alpha_1, \tilde{\alpha}_1}^{\mu}$, which we assume to be normalized to $\int f \,d\mu = 1$ without loss of generality, and we show that $L_{K, \mu} f \in \mathcal{F}^{\nu}_{\alpha_2, \tilde{\alpha}_2}$.
		For $x \in B_{m}$, we find that $L_{K, \mu}(f)(x) \geq 0$ holds using Equation \eqref{eq:gengeqclosed}, since
		\begin{align*}
		L_{K, \mu}f(x) &= \int K(x, y) f(y) \,\mu(dy) \\
		&\geq \int K(x, y) f_+(y) \,\mu(dy) - \int f_- \,d\mu \\
		&\geq l(m) \int_{B_{m}} f_+(y) \,\mu(dy) - \tilde{\alpha}_1(m) \\
		&\geq l(m)(1-\alpha_1(m)) - \tilde{\alpha}_1(m) \geq 0.
		\end{align*}
		
		Next, we note that, for all $r \geq m$, we have $\int_{B_r} f_+ \,d\mu \geq (1-\alpha_1(r))$ and thus with a similar calculation to the one above
		\begin{align*}
		\int L_{K, \mu}f\,d\nu \geq  l(r) \nu(B_r)\int_{B_r} f_+ \,d\mu - \tilde{\alpha}_1(m) \geq l(r) \nu(B_r) (1-\alpha_1(r)) - \tilde{\alpha}_1(m)
		\end{align*}
		and taking the supremum over $r \geq m$ as in \eqref{eq:genintegralcomparable} thus 
		\[
		\int L_{K, \mu} f \,d\nu \geq \bar{a}.
		\]
		Using this inequality, we next show the defining inequalities for $L_{K, \mu} f \in \mathcal{F}^{\nu}_{\alpha_2, \tilde{\alpha}_2}$. For $r \geq m$, using Equation \eqref{eq:gentailplusclosed} and Jensen's inequality for the function $x \mapsto x_+$, we have
		\begin{align*}
		\int_{B_r^{\com}} (L_{K, \mu} f)_+\,d\nu &\leq \int_{B_r^{\com}} \int f_+(y) \,\mu(dy) \,\nu(dx)\\ 
		&\leq (1+\tilde{\alpha}_1(m)) \nu(B_r^{\com})\\
		&\leq \frac{1+\tilde{\alpha}_1(m)}{\bar{a}} \, \nu(B_r^{\com}) \int L_{K, \mu} f \,d\nu \\
		&\leq \alpha_2(r) \int L_{K, \mu} f \,d\nu.
		\end{align*}
		And similarly, for $r \geq m$, using Equation \eqref{eq:gentailminusclosed} and Jensen's inequality for the function $x \mapsto x_-$, 
		\begin{align*}
		\int_{B_r^{\com}} (L_{K, \mu} f)_- &\leq \int_{B_r^{\com}} \int f_-(y) \,\mu(dy) \,\nu(dx) \\
		&\leq \tilde{\alpha}_1(m) \nu(B_r^{\com}) \\
		&\leq \frac{\tilde{\alpha}_1(m)}{\bar{a}} \nu(B_r^{\com}) \int L_{K, \mu} f \,d\nu\\
		&\leq \tilde{\alpha}_2(r) \int L_{K, \mu} f \,d\nu.
		\end{align*}
		Thus, we have shown that $L_{K, \mu}(f) \in \mathcal{F}^{\nu}_{\alpha_2, \tilde{\alpha}_2}$.

		The second part of the statement, i.e., that $L_{K^\top, \nu}(\mathcal{G}_{\alpha_2, \tilde{\alpha}_2}^{\nu}) \subseteq \mathcal{G}_{\alpha_1, \tilde{\alpha}_1}^{\mu}$, follows simply from the first part. Indeed, for $f \in \mathcal{F}^{\mu}_{\alpha_1, \tilde{\alpha}_1}$ and $g \in \mathcal{G}^{\nu}_{\alpha_2, \tilde{\alpha}_2}$, we have
		\[
		\int L_{K^\top, \nu}(g) f \,d\mu = \int L_{K, \mu}(f) g \,d\nu \geq 0,
		\]
		where the final inequality holds since $L_{K, \mu}(f) \in \mathcal{F}^{\nu}_{\alpha_2, \tilde{\alpha}_2}$.
	\end{proof}
\end{proposition}

The above Proposition \ref{prop:nonexpansive} only yields that the map $L_{K, \mu}$ is non-expansive (i.e., a contraction with coefficient 1). In the following, we give a crude yet effective way to show that it is also a strict contraction (i.e., with contraction coefficient strictly lower than 1) under suitable assumptions. 

The following gives the basic approach used to bound the diameter $\Delta$ from Birkhoff's contraction theorem, cf.~Theorem \ref{thm:birkhoff}. 
\begin{lemma}\label{lem:approachbound}
	We have
	\begin{align*}
	\sup_{g, \tilde{g} \in \mathcal{G}_{\alpha_2, \tilde{\alpha}_2}^{\nu, norm}} d_{\mathcal{G}^{\mu}_{\alpha_1, \tilde{\alpha}_1}}\left(L_{K^{\top}, \nu}(g), L_{K^{\top}, \nu}(\tilde{g})\right)
	\\\leq 2 \Bigg( \log \sup_{f\in \mathcal{F}_{\alpha_1, \tilde{\alpha}_1}^{\mu, norm}} \sup_{g \in \mathcal{G}_{\alpha_2, \tilde{\alpha}_2}^{\nu, norm}} \int L_{K^{\top}, \nu}(g) f \,d\mu \\-  \log \inf_{f\in \mathcal{F}_{\alpha_1, \tilde{\alpha}_1}^{\mu, norm}} \inf_{g \in \mathcal{G}_{\alpha_2, \tilde{\alpha}_2}^{\nu, norm}} \int L_{K^{\top}, \nu}(g) f \,d\mu \Bigg).
	\end{align*}
	\begin{proof}
		The statement will follow from Lemma \ref{lem:calcdist}. To this end, we note that the supremum therein may without loss of generality be taken over $f \in  \mathcal{F}_{\alpha_1, \tilde{\alpha}_1}^{\mu, norm}$, since $\int f \,d\mu = 0$ would imply $f = 0$ ($\mu$-a.s.). Thus, Lemma \ref{lem:calcdist} yields the statement since for all $g, \tilde{g} \in \mathcal{G}_{\alpha_2, \tilde{\alpha}_2}^{\nu, norm}$, we have
		\[
		\sup_{f \in \mathcal{F}_{\alpha_1, \tilde{\alpha}_1}^{\mu, norm}} \frac{\int L_{K^{\top}, \nu}(g) f\,d\mu}{\int L_{K^{\top}, \nu}(\tilde{g}) f \,d\mu} \leq \frac{\sup_{f \in \mathcal{F}_{\alpha_1, \tilde{\alpha}_1}^{\mu, norm}} \int L_{K^{\top}, \nu}(g) f\,d\mu}{\inf_{f \in \mathcal{F}_{\alpha_1, \tilde{\alpha}_1}^{\mu, norm}} \int L_{K^{\top}, \nu}(\tilde{g}) f\,d\mu}.
		\]
	\end{proof}
\end{lemma}

\begin{lemma}\label{lem:suffbounds}
	If $\int g \,d\nu > 0$ for all $g \in \mathcal{G}_{\alpha_2, \tilde{\alpha}_2}^{\nu} \setminus\{0\}$ and
	\begin{align*}
	\inf_{f\in \mathcal{F}_{\alpha_1, \tilde{\alpha}_1}^{\mu, norm}} \inf_{g \in \mathcal{G}_{\alpha_2, \tilde{\alpha}_2}^{\nu, norm}} \int L_{K^{\top}, \nu}(g) f \,d\mu &> 0,\\ \sup_{f\in \mathcal{F}_{\alpha_1, \tilde{\alpha}_1}^{\mu, norm}} \sup_{g \in \mathcal{G}_{\alpha_2, \tilde{\alpha}_2}^{\nu, norm}} \int L_{K^{\top}, \nu}(g) f \,d\mu &< \infty,
	\end{align*}
	then $L_{K^{\top}, \nu}$ is a strict contraction from $(\mathcal{G}^{\nu}_{\alpha_2, \tilde{\alpha}_2} \setminus\{0\}, d_{\mathcal{G}^{\nu}_{\alpha_2, \tilde{\alpha}_2}})$ into $(\mathcal{G}^{\mu}_{\alpha_1, \tilde{\alpha}_1}, d_{\mathcal{G}^{\mu}_{\alpha_1, \tilde{\alpha}_1}})$.
	\begin{proof}
		We show that $\Delta$ from Theorem \ref{thm:birkhoff} for $C_1 = \mathcal{G}^{\nu}_{\alpha_2, \tilde{\alpha}_2}$ and $C_2 = \mathcal{G}^{\mu}_{\alpha_1, \tilde{\alpha}_1}$ is finite. Indeed, since $\int g \,d\nu > 0$ for all $g \in \mathcal{G}_{\alpha_2, \tilde{\alpha}_2}^{\nu}\setminus\{0\}$, we can without loss of generality restrict to normalized functions $g \in \mathcal{G}_{\alpha_2, \tilde{\alpha}_2}^{\nu, norm}$ to calculate $\Delta$. Thus, applying Lemma \ref{lem:approachbound}, together with the remaining assumptions of the current Lemma, yields the claim.
	\end{proof}
\end{lemma}

The previous results given by Lemma \ref{lem:approachbound} and Lemma \ref{lem:suffbounds} yield abstract sufficient conditions for $L_{K^{\top}, \nu}$ to be a strict contraction. In the following, we show how to derive simple sufficient conditions depending only on the tail behavior of the marginals and the decay of the kernel. 
\begin{lemma}\label{lem:lowerandupper}
	For $f \in \mathcal{F}_{\alpha_1, \tilde{\alpha}_1}^{\mu, norm}$ and $g \in \mathcal{G}_{\alpha_2, \tilde{\alpha}_2}^{\nu, norm}$ and $r > 0$, we have 
	\begin{align*}
	\int L_{K^\top, \nu}(g) f \,d\mu \geq ~&\frac{1+l(r)}{2} \\
	&- \frac{1-l(r)}{2} \int_{B_r} |f| \,d\mu \int_{B_r} |g| \,d\nu \\
	&- \int |f| \,d\mu \int_{B_r^{\com}} |g|\,d\nu  \\
	&- \int_{B_r^{\com}} |f| \,d\mu \int_{B_r} |g|\,d\nu
	\end{align*}
	and
	\[
	\int L_{K^\top, \nu}(g) f \,d\mu \leq \int |f| \,d\mu \int |g|\,d\nu.
	\]
	\begin{proof}
		The second inequality follows from the normalization $\|K\|_\infty \leq 1$.
		
		For the first inequality, let $\bar{e} := (1+l(r))/2$. We have 
		\begin{align*}
		&\iint_{B_r \times B_r} (K(x, y) - \bar{e}) f(x) g(y) \,\nu(dy)\mu(dx)\\ &\geq - \sup_{x, y \in B_r} |K(x, y) - \bar{e}| \, \int_{B_r} |f| \,d\mu \int_{B_r} |g|\,d\nu \\ &\geq -\frac{1-l(r)}{2} \, \int_{B_r} |f| \,d\mu \int_{B_r} |g|\,d\nu.
		\end{align*}
		Further,
		\begin{align*}
		&\iint_{(B_r \times B_r)^{\com}} (K(x, y) - \bar{e}) f(x) g(y) \,\nu(dy)\mu(dx)\\
		&\geq - \iint_{(B_r \times B_r)^{\com}} |f(x)| |g(y)| \,\nu(dy) \mu(dx) \\
		&\geq -\int |f| \,d\mu \int_{B_r^{\com}} |g|\,d\nu- \int_{B_r^{\com}} |f| \,d\mu \int_{B_r} |g|\,d\nu.
		\end{align*}
		We conclude by noting
		\[
		\int L_{K^\top, \nu}(g) f \,d\mu = \bar{e} + \iint (K(x, y) - \bar{e}) f(x)g(y) \,\nu(dy)\mu(dx),
		\]
		splitting the latter integral into $(B_r\times B_r) \stackrel{.}{\cup} (B_r \times B_r)^{\com}$ and plugging in the above obtained lower bounds.
	\end{proof}
\end{lemma}

We emphasize that the above Lemma \ref{lem:lowerandupper} is a rather crude way to obtain the lower bound of interest (which is the main factor determining the contraction coefficient). The approach is crude since it basically boils down to approximating the kernel $K$ by a constant $\bar{e}$ on a large compact $B_r \times B_r$, and to cut off the remaining tails. It is hence almost surprising that this crude approach still leads to such strong results like Theorem \ref{thm:generalkernel} and Theorem \ref{thm:sinkhorn}. We at least believe that the resulting constants---which we do not track explicitly here---could be significantly improved if a more fine-grained approach could be found. We leave such potential refinements for future work.

The following uses the notation $\beta_2$ and $\tilde{\beta}_2$ for $g \in \mathcal{G}_{\alpha_2, \tilde{\alpha}_2}^{\nu}$ as defined prior to Lemma \ref{lem:tailboundsfg}, and makes the bound from the above Lemma \ref{lem:lowerandupper} more explicit.
\begin{lemma}\label{lem:concretelowerbound}
	For $f \in \mathcal{F}_{\alpha_1, \tilde{\alpha}_1}^{\mu, norm}$ and $g \in \mathcal{G}_{\alpha_2, \tilde{\alpha}_2}^{\nu, norm}$ and $r\geq m$, assuming that $\tilde{\alpha}_1(m) \leq 0.5$, $\tilde{\beta}_2(m) \leq 0.5$ and $l(r) \leq 1$, we get
	\[
	\int L_{K^\top, \nu}(g) f \,d\mu \geq l(r) - \tilde{\alpha}_1(m) - 4 \tilde{\beta}_2(m) - 4 \alpha_1(r) - 4 \tilde{\alpha}_1(r) - 5 \beta_2(r) - 5 \tilde{\beta}_2(r).
	\]
	\begin{proof}
		We can first rearrange the bounds from Lemma \ref{lem:lowerandupper} to get
		\begin{align*}
		\int L_{K^\top, \nu}(g) f \,d\mu \geq &~\frac{l(r)}{2} \left(1+\int_{B_r} |f|\,d\mu \int_{B_r} |g| \,d\nu\right) \\
		&- \frac{1}{2} \left(\int_{B_r} |f|\,d\mu \int_{B_r} |g|\,d\nu - 1\right)\\
		&- \int |f| \,d\mu \int_{B_r^{\com}} |g|\,d\nu  \\
		&- \int_{B_r^{\com}} |f| \,d\mu \int_{B_r} |g|\,d\nu.
		\end{align*}
		Further, Lemma \ref{lem:tailboundsfg} yields (using, at various points, that $\frac{1}{1-\tilde{\beta}_2(m)} \leq 2$)
		\begin{align*}
		\int_{B_r} |f|\,d\mu \int_{B_r} |g|\,d\nu &\leq (1+2\tilde{\alpha}_1(m)) \Big(1+2\frac{\tilde{\beta}_2(m)}{1-\tilde{\beta}_2(m)}\Big) \leq 1 + 2 \tilde{\alpha}_1(m) + 8 \tilde{\beta}_2(m), \\
		\int |f|\,d\mu \int_{B_r^{\com}} |g|\,d\nu &\leq (1+2\tilde{\alpha}_1(m))\Big(\frac{\tilde{\beta}_2(r) + \beta_2(r)}{1-\tilde{\beta}_2(m)}\Big) \leq 4 \left(\tilde{\beta}_2(r) + \beta_2(r)\right),\\
		\int_{B_r^{\com}} |f|\,d\mu \int_{B_r} |g|\,d\nu &\leq (\alpha_1(r) + \tilde{\alpha}_1(r)) \Big(1+2\frac{\tilde{\beta}_2(m)}{1-\tilde{\beta}_2(m)}\Big) \leq 3 \left(\alpha_1(r) + \tilde{\alpha}_1(r)\right), \\
		\int_{B_r} |f|\,d\mu \int_{B_r} |g| \,d\nu &\geq 1 - \alpha_1(r) - \tilde{\alpha}_1(r) - 2\beta_2(r) - 2\tilde{\beta}_2(r).
		\end{align*}
		Plugging these into the first inequality of the proof  and using $l(r) \leq 1$ yields the claim.
	\end{proof}
\end{lemma}

Finally, we are ready to state the main result of this section, which gives conditions for $L_{K^{\top}, \nu}$ to be a contraction from $\mathcal{G}^\nu_{\alpha_2, \tilde{\alpha}_2}$ into $\mathcal{G}^{\mu}_{\alpha_1, \tilde{\alpha}_1}$. The conditions only depend on the tail functions of the respective cones, the tails of $\nu$, and the lower decay $l$ of $K$. Regarding notation, we recall that for two functions $f_1, f_2 : (0, \infty) \rightarrow [0, \infty)$, we write $f_1 \in o(f_2)$ if $\lim_{r \rightarrow \infty} \frac{f_1(r)}{f_2(r)} = 0$ (with the usual convention $0/0 = \infty$). \SE{We further write $\nu(B_\cdot^\com)$ for the function $(0, \infty) \ni r \mapsto \nu(B_r^\com)$.}
\begin{theorem}\label{thm:generalkernel}
	Assume \[
	\alpha_1, \tilde{\alpha}_1, \frac{\nu(B_\cdot^\com)}{\alpha_2}, \frac{\nu(B_\cdot^\com)}{\tilde{\alpha}_2} \in o(l).\]
	Then, there exists $\underline{m} > 0$ such that for all $m \geq \underline{m}$, there exists $\kappa \in (0, 1)$ such that for all $g, \tilde{g} \in \mathcal{G}_{\alpha_2, \tilde{\alpha}_2}^{\nu}\setminus\{0\}$, we have
	\[
	d_{\mathcal{G}_{\alpha_1, \tilde{\alpha}_1}^{\mu}}\left( L_{K^{\top}, \nu}(g), L_{K^{\top}, \nu}(\tilde{g}) \right) \leq \kappa d_{\mathcal{G}_{\alpha_2, \tilde{\alpha}_2}^{\nu}}\left(g, \tilde{g}\right).
	\]
	\begin{proof}
		We show that the assumptions of Lemma \ref{lem:suffbounds} are satisfied by using the bounds from Lemma \ref{lem:concretelowerbound}. The upper bound in Lemma \ref{lem:suffbounds} is of course always satisfied \SE{since an upper bound of $1$ follows using $\|K\|_\infty \leq 1$}. The only thing to prove is the lower bound in Lemma \ref{lem:suffbounds}, and that $\int g \,d\mu > 0$ for all $g \in \mathcal{G}^{\nu}_{\alpha_2, \tilde{\alpha}_2}\setminus\{0\}$. 
		
		First, we prove $\int g \,d\mu > 0$ for all $g \in \mathcal{G}^{\nu}_{\alpha_2, \tilde{\alpha}_2}$ for large enough $\underline{m}$. Take $\beta_2$ and $\tilde{\beta}_2$ as defined prior to Lemma \ref{lem:tailboundsfg} and assume $\underline{m}$ is large enough so that $\tilde{\beta}_2(\underline{m}) < 1$. Then, Lemma \ref{lem:tailboundsfg} (ii) yields
		\[
		\int_{B_m^{\com}} |g|\,d\mu \leq \frac{\beta_2(m) + \tilde{\beta}_2(m)}{1-\tilde{\beta}_2(m)} \int g \,d\mu.
		\]
		Hence, if $\int g \,d\mu = 0$ holds, the above inequality, together with $g \eins_{B_m} \geq 0$, implies $g = 0$ $\mu$-a.s..
		
		We now turn to the lower bound in Lemma \ref{lem:suffbounds}.
		Let $\varepsilon < \frac{1}{27}$ and take $\underline{m}$ large enough such that $\tilde{\alpha}_2(\underline{m}) \leq 0.5$ and, for all $r \geq \underline{m}$, we have
		\begin{align*}
		\max\left\{\alpha_1(r) , \tilde{\alpha}_1(r) , \frac{t_{\nu}(r)}{\alpha_2(r)} , \frac{t_{\nu}(r)}{\tilde{\alpha}_2(r)}\right\} &\leq \varepsilon l(r) \\
		\frac{1}{\nu(B_{\underline{m}})} &\leq 1+\varepsilon.
		\end{align*}
		Let now $m \geq \underline{m}$.
		The above implies, by definition of $\beta_2$ and $\tilde{\beta}_2$, that
		\begin{align*}
		\tilde{\beta}_2(r) &= \frac{1}{\nu(B_{r})} \sup_{a \geq r} \frac{\nu(B_a^\com)}{\alpha_2(a)} \leq (1+\varepsilon) \varepsilon l(r), ~~\text{for all }r \geq m,\\
		\beta_2(r) &= \frac{1}{\nu(B_r)} \Big(\sup_{a \geq r} \frac{\nu(B_a^\com)}{\tilde{\alpha}_2(a)} + \nu(B_r^\com)\Big) \leq \frac{3}{2}(1+\varepsilon)  \varepsilon l(r), ~~\text{for all }r \geq m,
		\end{align*}
		where the last inequality uses $\tilde{\alpha}_2(r) \leq 0.5$ (and thus $\nu(B_r^\com) \leq 0.5 \varepsilon l(r)$).
		
		Plugging in the bound of Lemma \ref{lem:concretelowerbound} for $f \in \mathcal{F}_{\alpha_1, \tilde{\alpha}_1}^{\mu, norm}$ and $g \in \mathcal{G}_{\alpha_2, \tilde{\alpha}_2}^{\nu, norm}$ and $r = m$, we obtain
		\[\int L_{K^\top, \nu}(g) f \,d\mu \geq l(m) (1 - 9 \varepsilon - 16.5 (1+\varepsilon) \varepsilon) > 0,
		\]
		where the final inequality holds since $\varepsilon < \frac{1}{27}$. We have thus shown that Lemma \ref{lem:suffbounds} can be applied, which yields the claim.
	\end{proof}
\end{theorem}

While Theorem \ref{thm:generalkernel} is stated generally and with simple conditions, the given result may still seem slightly abstract because the tail functions are regarded as exogenously given, while they should more accurately be seen as a modelling tool which can be specified depending on the other parameters of the setting. Hence, the following corollary shows how one may choose the tail functions in the important case $\mu=\nu$.
\begin{corollary}\label{cor:cleankernel}
	Assume $\lim_{r \rightarrow \infty} \frac{\mu(B_r^{\com})}{l(r)^2} = 0.$ Then, there exists $\underline{m} > 0$ such that for all $m > \underline{m}$ and $\alpha(r) = \tilde{\alpha}(r) = \sqrt{\mu(B_r^{\com})}$, $r\geq 0$, we have that $L_{K, \mu}$ is a strict contraction on $\mathcal{G}_{\alpha, \tilde{\alpha}}^{\mu}\setminus\{0\}$ with respect to Hilbert's metric corresponding to the cone $\mathcal{G}_{\alpha, \tilde{\alpha}}^{\mu}$.
	\begin{proof}
		Note that using $K^{\top}$ or $K$ doesn't matter because the assumption on $K$ through $l$ is symmetric. The claim thus follows from Theorem \ref{thm:generalkernel} with $\mu = \nu$ and $\alpha_1=\alpha_2=\tilde{\alpha}_1=\tilde{\alpha}_2=\alpha$.
	\end{proof}
\end{corollary}

\section{Sinkhorn's algorithm}\label{sec:Sinkhorn}

We are given two marginal distributions $\mu, \nu \in \mathcal{P}(\X)$.\footnote{Using different spaces $\X_1$ and $\X_2$ works analogously, but we aim to keep notation simple. In fact, using the same space for both marginals can be done without loss of generality, since the case in which $\tilde{\mu} \in \mathcal{P}(\X_1)$ and $\tilde{\nu} \in \mathcal{P}(\X_2)$ can be derived from the given setting by using $\X = \X_1 \times \X_2$, $\mu = \tilde{\mu} \otimes \delta_{x_{0, 1}} \in \mathcal{P}(\X)$ and $\nu = \delta_{x_{0, 2}} \otimes \tilde{\nu} \in \mathcal{P}(\X)$ for suitable points $x_{0, i} \in \X_i$, $i=1, 2$.}
Sinkhorn's algorithm is a concatenation of several simple operations. Starting with $g^{(0)} : \X \rightarrow \mathbb{R}_+$ (e.g., $g^{(0)} \equiv 1$) the steps of Sinkhorn's algorithm are given by
\begin{align}\label{eq:Sinkhornstep}
g^{(n+1)} = \left(S \circ L_{K^{\top}, \nu} \circ S \circ L_{K, \mu}\right) \big(g^{(n)}\big) ~~ \text{ for } n \in \mathbb{N},
\end{align}
where $S(g) := 1/g$. For a motivation of Sinkhorn's algorithm, we refer to Section \ref{subsec:intromain}.
We emphasize that the way we state Sinkhorn's algorithm is in terms of the $\mu$ marginal and multiplicative dual functions, meaning that $g^{(n)}$ should converge to $g_1^*$, where $g_1^*$ is an optimal multiplicative dual potential for the $\mu$ marginal. Further, the kernel $K(x, y) = \exp(-c(x, y)/\varepsilon)$ is given by the cost function. 

While kernel integral operators were already treated in the previous section, for a thorough analysis of Sinkhorn's algorithm we need to further understand the operation $S(g) = 1/g$. Ideally, we want to show that $S$ is a contraction, or at least non-expansive.
To this end, the main idea is to use the elementary identity
\[
d_{\mathcal{G}^{\mu}_{\alpha, \tilde{\alpha}}}\left( S(g), S(\tilde{g}) \right) = d_{\mathcal{G}^{\mu}_{\alpha, \tilde{\alpha}}}\left( \frac{\tilde{g}}{g \tilde{g}}, \frac{g}{g\tilde{g}}\right).
\]
This identity reveals that understanding the operation $S$ will be a special case of understanding how Hilbert's metric between two functions $g$ and $\tilde{g}$ changes when both functions are multiplied by some function $h$ (in the case of $S$, we have $h = \frac{1}{g \tilde{g}}$). In this regard, it will follow from Lemma \ref{lem:multclosednew} (together with Lemma \ref{lem:calcdist}) that
\[
d_{\mathcal{G}^{\mu}_{\alpha_1, \tilde{\alpha}_1}}\left(g h, \tilde{g} h \right) \leq d_{\mathcal{G}^{\mu}_{\alpha_2, \tilde{\alpha}_2}}(g, \tilde{g})
\]
for a suitable relation between $\alpha_1, \tilde{\alpha}_1$, $h$ and $\alpha_2, \tilde{\alpha}_2$.
While this will not quite yield that $S$ is a contraction, it will at least yield that it is non-expansive when the right hand side uses a (slightly) larger set of test functions.
To apply Lemma \ref{lem:multclosednew} under suitable assumptions on $h$, Lemma \ref{lem:growthnew} shows the growth behavior of functions resulting from Sinkhorn's algorithm. Lemma \ref{lem:multclosedexp} simplifies Lemma \ref{lem:multclosednew} in a setting with exponential growth. Proposition \ref{prop:inversion} shows that the operator $S$ is non-expansive up to a change of the set of test functions.  The main result of this section, Theorem \ref{thm:sinkhorn}, shows contractivity of Sinkhorn's iterations with respect to Hilbert's metric, which yields exponential convergence as an immediate consequence. Finally, Corollary \ref{cor:Sinkhorn} transfers the convergence from Hilbert's metric to the total variation norm of the primal iterations, which builds on Proposition \ref{prop:relationmetrics} and an embedding into the product space given in Lemma \ref{lem:embednew}.

In the following, for a non-decreasing function $\kappa: \mathbb{R}_+ \rightarrow \mathbb{R}_+$, we use its generalized inverse $\kappa^{-1}(r) := \inf\{a \geq 0 : r \leq \kappa(a)\}$, $r \geq 0$.\footnote{Leading to the fact that $\kappa^{-1}(r) \leq a \Leftrightarrow r \leq \kappa(a)$ for $a, r \geq 0$.}

\begin{lemma}\label{lem:multclosednew}
	Let $f \in \mathcal{F}^{\mu}_{\alpha, \tilde{\alpha}}$ and $h: \X \rightarrow \mathbb{R}_+$. Then, we have $fh \in \mathcal{F}^{\mu}_{\gamma, \tilde{\gamma}}$ in both of the following cases.
	\begin{itemize}
		\item[(i)] Let $\kappa : \mathbb{R}_+ \rightarrow \mathbb{R}_+$ be non-decreasing and unbounded, and assume $1 \leq h(x) \leq \kappa(d_{\X}(x_0, x))$ for all $x \in \X$. Let $\alpha, \tilde{\alpha}$ be continuous and define
		\begin{align*}
		\xi(r) := \int_0^\infty \alpha(\kappa^{-1}(s)\lor r)\,ds ~&\text{ and } ~ \tilde{\xi}(r) := \int_0^\infty \tilde{\alpha}(\kappa^{-1}(s)\lor r) \,ds,\\
		\gamma(r) = \frac{\xi(r)}{1-\tilde{\xi}({m})} ~&\text{ and }~
		\tilde{\gamma}(r) = \frac{\tilde{\xi}(r)}{1-\tilde{\xi}({m})}, ~~ r \geq 0,
		\end{align*}
		where we assume $\tilde{\xi}(m) < 1$.
		\item[(ii)] Let $l : \mathbb{R}_+ \rightarrow [0, 1)$ be non-increasing, assume $l(d_{\X}(x_0, x)) \leq h(x) \leq 1$ for all $x \in \cX$, $\alpha(m) < 1$ and $\left(1- \frac{\tilde{\alpha}(m)}{l(m) (1-\alpha(m))}\right) > 0$, and define
		\begin{align*}
		\gamma(r) &= \frac{\alpha(r)}{\zeta(m)} ~\text{ and }~ \tilde{\gamma}(r) = \frac{\tilde{\alpha}(r)}{\zeta(m)}, \text{ where}\\
		\zeta(m) &= l(m) (1-\alpha(m)) \left(1- \frac{\tilde{\alpha}(m)}{l(m) (1-\alpha(m))}\right).
		\end{align*}
	\end{itemize}
	\begin{proof}
		We start with case (i): Let $r \geq m$. Then, using the layer cake representation, we have
		\begin{align}\label{eq:fmc1}
		\begin{split}
		\int_{B_r^{\com}} (fh)_+ \,d\mu &= \int_{B_r^{\com}} f_+ h \,d\mu \\
		&\leq \int_{B_r^{\com}} f_+(x) \kappa(d_{\X}(x_0, x)) \,\mu(dx) \\
		&= \int_{B_r^{\com}} \int_{0}^\infty \eins_{\{s \leq \kappa(d_{\X}(x_0, x))\}} \,ds f_+(x) \,\mu(dx) \\
		&= \int_0^\infty \int_{B_r^{\com} \cap \{d_{\X}(x, x_0) \geq \kappa^{-1}(s)\}} f_+(x) \,\mu(dx) \,ds \\
		&\leq \int_0^\infty \alpha(\kappa^{-1}(s) \lor r) \,ds \int f \,d\mu \\
		&=\xi(r) \int f \,d\mu,
		\end{split}
		\end{align}
		where we note that the second inequality uses continuity of $\alpha$ (to allow for the fact that we have $d_{\X}(x, x_0) \geq \kappa^{-1}(s)$ instead of $d_{\X}(x, x_0) > \kappa^{-1}(s)$).
		For $r\geq m$, we analogously obtain 
		\begin{align}\label{eq:fmc2}
		\int_{B_r^{\com}} (fh)_- \,d\mu \leq \tilde{\xi}(r) \int f \,d\mu.
		\end{align}
		From the latter and $h \geq 1$ we further find
		\begin{equation}\label{eq:fmc3}
		\int fh \,d\mu \geq \int f_+ h \,d\mu - \tilde{\xi}({m}) \int f \,d\mu \geq (1-\tilde{\xi}({m})) \int f \,d\mu.
		\end{equation}
		Concatenating each of \eqref{eq:fmc1} and \eqref{eq:fmc2} with \eqref{eq:fmc3} completes the proof of (i).
		
		Next, we turn to case (ii). We have that 
		\[
		\int_{B_r^{\com}} (fh)_+ \,d\mu \leq \alpha(r) \int f \,d\mu ~ \text{ and } ~ \int_{B_r^{\com}} (fh)_- \,d\mu \leq \tilde{\alpha}(r) \int f \,d\mu.
		\]
		Hence, it only remains to show that $\int f \,d\mu \leq \frac{1}{\zeta(m)} \int fh \,d\mu$. To this end, note first that we have
		\begin{align*}
		\int f \,d\mu &\leq \frac{1}{(1-\alpha(m))} \int_{B_{m}} f_+ \,d\mu,\\
		\int_{B_{m}} f_+ \,d\mu &\leq \frac{1}{l(m)} \int_{B_{m}} (fh)_+ \,d\mu.
		\end{align*}
		Using these two inequalities, we obtain
		\begin{align*}
		\int fh \,d\mu &= \int (fh)_+ \,d\mu - \int_{B_{m}^{\com}} (fh)_- \,d\mu \\
		&\geq  \int (fh)_+ \,d\mu - \tilde{\alpha}(m) \int f\,d\mu \\
		&\geq \int (fh)_+ \,d\mu - \frac{\tilde{\alpha}(m)}{l(m) (1-\alpha(m))} \int_{B_{m}} (fh)_+\,d\mu \\
		&\geq \left(1-\frac{\tilde{\alpha}(m)}{l(m) (1-\alpha(m))}\right) \int_{B_{m}} (fh)_+ \,d\mu \\
		&\geq l(m) \left(1-\frac{\tilde{\alpha}(m)}{l(m) (1-\alpha(m))}\right) \int_{B_{m}} f_+ \,d\mu \\
		&\geq \zeta(m) \int f \,d\mu,
		\end{align*}
		which yields the claim.
	\end{proof}
\end{lemma}

The following result is the main tool to bound the growth of the functions $h$ for which Lemma \ref{lem:multclosednew} will be applied later on. To this end, we define
\begin{equation}\label{eq:defcalA}
\mathcal{A} := \{L_{K, \mu}(g) : g \in \mathcal{G}_{\alpha_1, \tilde{\alpha}_1}^{\mu}, g \geq 0\},
\end{equation}
which can be seen as the set of functions which may occur during Sinkhorn's algorithm as an input to the operator $S$ (along with the obvious analogue defined through $L_{K^\top, \nu}$).
\begin{lemma}\label{lem:growthnew} Define $\iota_1 : \mathbb{R}_+ \rightarrow \mathbb{R}_+$ by
	\begin{align*}
	\iota_1(r) &:= \frac{l(r\lor m)}{1+\beta_1(m)}
	\end{align*}
	Then, we have $\mathcal{A} \subseteq \{\lambda g : \lambda > 0, \iota_1(r) \leq g(x) \leq 1 \text{~for all~} r\geq0 \text{ and } x \in B_r\}$.
	\begin{proof}
		Let $L_{K, \mu}(g) \in \mathcal{A}$ for some $g \in \mathcal{G}_{\alpha_1, \tilde{\alpha}_1}^{\mu}$, $g \geq 0$. Clearly, since $\|K\|_\infty \leq 1$, this implies $L_{K, \mu}(g) \leq \int g \,d\mu$. 
		
		We note that Lemma \ref{lem:tailboundsfg} (i) implies, for all $r \geq m$,
		\[
		\int_{B_r} g \,d\mu 
		\geq \int g\,d\mu - \beta_1(r) \int_{B_{m}} g\,d\mu \geq \int g \,d\mu - \beta_1(r) \int_{B_r} g \,d\mu
		\]
		and thus $(1+\beta_1(m)) \int_{B_r} g\,d\mu \geq \int g \,d\mu$ (note that $\beta_1(r) \leq \beta_1(m)$). Hence, we get, for any $x \in B_r$,
		\begin{align*}
		L_{K, \mu}(g)(x) &= \int K(y, x) g(y) \,\mu(dy) \\
		&\geq \int_{B_{r\lor m}} K(y, x) g(y) \,\mu(dy) \\
		&\geq l(r \lor m) \int_{B_{r\lor m}} g \,d\mu \\
		&\geq \frac{l(r \lor m)}{1+\beta_1(m)} \int g\,d\mu.
		\end{align*}
		Together, we find that with $\lambda := \int g \,d\mu$ and $\tilde{g} := L_{K, \mu}(g)/\lambda$, we have that $L_{K, \mu}(g) = \lambda \tilde{g}$ and $\iota_1(r) \leq \tilde{g}(x) \leq 1$ for all $x \in B_r, r \geq 0$. This completes the proof.
	\end{proof}
\end{lemma}

The following simplifies the application of Lemma \ref{lem:multclosednew} in a setting where all growth functions are exponential of a suitable order.
\begin{lemma}\label{lem:multclosedexp}
	Let $p_1 > p_2 \geq p_3 > 0$, $\iota(r) = \exp(-r^{p_3})$, $\alpha_i(r) = \exp(-r^{p_i})$ for $i=1, 2$ and $r \geq 0$. Then, there exists $\underline{m} > 0$ such that for all $m > \underline{m}$, we have $fh \in \mathcal{F}_{\alpha_2, \alpha_2}^{\mu}$ for all $f \in \mathcal{F}_{\alpha_1, \alpha_1}^{\mu}$ and $h : \cX \rightarrow \mathbb{R}$ satisfying either $1\leq h(x) \leq \frac{1}{\iota(d_{\X}(x_0, x))}$ or $\iota(d_{\X}(x_0, x)) \leq h(x) \leq 1$ for all $x \in B_m^{\com}$.
	\begin{proof}
		The given assumptions on $h$ become weaker for increasing $p_3$ and are hence weakest for $p_3 = p_2$ and thus we prove this case without loss of generality (we will thus replace each occurrence of $p_3$ by $p_2$).
		
		First, consider the case where $1 \leq h(x) \leq \kappa(d_{\X}(x_0, x))$ with $\kappa(r) := 1/\iota(r)$. We want to apply Lemma \ref{lem:multclosednew} (i).
		The inverse of $\kappa$ is given by $\kappa^{-1}(r) = \log(r)^{1/p_2}$. Further, we note that for $r \geq 0$ large enough,
		\[
		e^{r^{p_2}} e^{-r^{p_1}} \leq 0.4 e^{-r^{p_2}},
		\]
		since $p_1 > p_2$, which can be seen by writing $e^{r^{p_2}-r^{p_1}} = e^{r^{p_1}(r^{p_2}/r^{p_1} - 1)}$. Similarly, we get that
		\[
		\alpha_1(\kappa^{-1}(r)) = \exp(-\log(r)^{p_1/p_2}) \leq 0.4 \exp(-2 \log(r)) = 0.4 r^{-2}
		\]
		for $r$ large enough, since $p_1 > p_2$.
		Thus, we can bound $\xi$ and $\tilde{\xi}$ from Lemma \ref{lem:multclosednew} (i) for $r \geq m$ by
		\begin{align*}
		\xi(r) = \tilde{\xi}(r) &= \int_0^\infty \alpha_1(\kappa^{-1}(s) \lor r) \,ds \\
		&= \kappa(r) \alpha_1(r) + \int_{\kappa(r)}^{\infty} \alpha_1(\kappa^{-1}(s))\,ds \\
		&\leq 0.4 \exp(-r^{p_2}) + 0.4 (\kappa(r))^{-1} = 0.8 \exp(-r^{p_2})
		\end{align*}
		and hence
		\[
		\gamma(r) = \tilde{\gamma}(r) = \frac{\xi(r)}{1-\tilde{\xi}(m)} \leq \exp(-r^{p_2}) = \alpha_2(r)
		\]
		for $r \geq m$ and $m$ large enough. This completes the first part of the proof.
		
		Next, we treat the case where $\iota(d_{\X}(x_0, x)) \leq h(x) \leq 1$ by applying Lemma \ref{lem:multclosednew} (ii). First, note first that under the given assumptions, it is easy to see that for $m$ large enough, $\zeta(m) \geq 0.5 \iota(m)$. Thus,
		\[
		\gamma(r) = \tilde{\gamma}(r) \leq 2 \frac{\alpha_1(r)}{\iota(m)} \leq 2 \exp(-r^{p_1} + r^{p_2}) \leq \exp(-r^{p_2}) = \alpha_2(r)
		\]
		for $r$ large enough, similarly to the above using $p_1 > p_2$. This completes the proof.
	\end{proof}
\end{lemma}

Next, we can obtain the non-expansiveness result of the operator $S$. As mentioned, the drawback here is that we must slightly increase the set of test functions for the right hand side.
\begin{proposition}\label{prop:inversion}
	Assume there exists $p > 0$ and $\delta > 0$ such that 
	\[
	\lim_{r \rightarrow \infty} \frac{\sup_{x, y \in B_r} c(x, y)}{r^p} = 0 ~~~\text{and}~~~\lim_{r \rightarrow \infty} \frac{\max\{\mu(B_r^{\com}), \nu(B_r^{\com})\}}{\exp(-r^{p+\delta})} = 0.
	\]
	For $p < p_2 < p_1 < p + \delta$, let $\alpha_i(r) = \exp(-r^{p_i})$ for $r \geq 0$. Then, there exists $\underline{m} > 0$ such that for all $m \geq \underline{m}$, we have
	\[
	d_{\mathcal{G}^{\mu}_{\alpha_1, \alpha_1}}(S(g), S(\tilde{g})) \leq d_{\mathcal{G}^{\mu}_{\alpha_2, \alpha_2}}(g, \tilde{g}) ~~ \text{ for all } g, \tilde{g} \in \mathcal{A}.
	\]
	\begin{proof}
		First note that the assumption on $c$ means that $l$ satisfies
		\begin{equation}\label{eq:growthlsink}
		\lim_{r \rightarrow \infty} \frac{\exp(-r^p)}{l(r)} = 0.
		\end{equation}
		Since
		\[
		d_{\mathcal{G}^{\mu}_{\alpha_1, \alpha_1}}(1/g, 1/\tilde{g}) = d_{\mathcal{G}^{\mu}_{\alpha_1, \alpha_1}}\left(\frac{\tilde{g}}{g \tilde{g}}, \frac{g}{g \tilde{g}}\right),
		\]
		our goal is to apply Lemma \ref{lem:multclosedexp} with $h = 1/(g \tilde{g})$. To this end, we find that for $\iota_1$ from Lemma \ref{lem:growthnew}, we have 
		\[
		\iota_1(r) \geq \frac{l(r \lor m)}{2}
		\]
		for $m$ large enough and thus, for a possibly re-scaled version of $h$ (i.e., we possibly multiply $h$ by a constant to normalize the lower bound to 1),
		\[
		1 \leq h(x) \leq \frac{1}{\iota_1(d_{\X}(x_0, x))^2} \leq 4 \exp(2 (d_{\X}(x_0, x)\lor m)^{p}) \leq \exp(d_{\X}(x_0, x)^{p_2})
		\]
		for $x \in B_m^{\com}$ and $m$ large enough. Thus, Lemma \ref{lem:multclosedexp} is applicable which yields $fh \in \mathcal{F}^{\mu}_{\alpha_2, \alpha_2}$ for $f \in \mathcal{F}^{\mu}_{\alpha_1, \alpha_1}$. Thus, an application of Lemma \ref{lem:calcdist} yields the claim.
	\end{proof}
\end{proposition}

We are ready to state the main result of this section, which shows that each step of Sinkhorn's algorithm is a contraction with respect to Hilbert's metric corresponding to a cone of the form $\mathcal{G}^{\mu}_{\alpha_1, \tilde{\alpha}_1}$.
As mentioned before, we state the convergence only in terms of the $\mu$ marginal. Since the problem is symmetric in $\mu$ and $\nu$, an analogous statement for the multiplicative dual potentials of $\nu$ holds as well.
\begin{theorem}\label{thm:sinkhorn}
	Assume there exists $p > 0$ and $\delta > 0$ such that 
	\[
	\lim_{r \rightarrow \infty} \frac{\sup_{x, y \in B_r} c(x, y)}{r^p} = 0 ~~~\text{and}~~~\lim_{r \rightarrow \infty} \frac{\max\{\mu(B_r^{\com}), \nu(B_r^{\com})\}}{\exp(-r^{p+\delta})} = 0.
	\]
	Then, for $\alpha(r) = \exp(-r^{p+\delta/2})$, there exists $\underline{m} > 0$ such that for all $m \geq \underline{m}$, each step of Sinkhorn's algorithm, given by \eqref{eq:Sinkhornstep} and initialized with $g^{(0)} \in \mathcal{G}_{\alpha, \alpha}^\mu\setminus \{0\}$ satisfying $g^{(0)} \geq 0$, is a strict contraction with respect to Hilbert's metric corresponding to $\mathcal{G}_{\alpha, \alpha}^\mu$.
	\begin{proof}
		We set $p_1 := p+\delta/2$, $\alpha_1(r) = \alpha(r)$ and $\alpha_2(r) := \exp(-r^{p_2})$ with $p_2 \in (p, p_1)$.
		By combining Proposition \ref{prop:inversion} and Theorem \ref{thm:generalkernel}, we find that for $m$ large enough, there exists $\kappa \in (0, 1)$ such that for all $g, \tilde{g} \in G^{\nu}_{\alpha_1, \tilde{\alpha}_1}$,
		\begin{align*}
		d_{\mathcal{G}^{\mu}_{\alpha_1, \alpha_1}}(S \circ L_{K, \mu}(g), S \circ L_{K, \mu}(\tilde{g})) &\leq d_{\mathcal{G}^{\mu}_{\alpha_2, \alpha_2}}(L_{K, \mu}(g), L_{K, \mu}(\tilde{g})) \\
		&\leq \kappa\, d_{\mathcal{G}^{\nu}_{\alpha_1, \alpha_1}}(g, \tilde{g}).
		\end{align*}
		And by symmetry, it analogously holds for all $g, \tilde{g} \in G^{\mu}_{\alpha_1, \tilde{\alpha}_1}$
		\begin{align*}
		d_{\mathcal{G}^{\nu}_{\alpha_1, \alpha_1}}(S \circ L_{K^{\top}, \nu}(g), S \circ L_{K^{\top}, \nu}(\tilde{g})) \leq \kappa\, d_{\mathcal{G}^{\mu}_{\alpha_1, \alpha_1}}(g, \tilde{g}).
		\end{align*}
		Combining both yields the claim.
	\end{proof}
\end{theorem}
We saw that the proof Theorem \ref{thm:sinkhorn} crucially uses that Theorem \ref{thm:generalkernel} can be applied with two different levels of growth for the test functions, thus reconciling the slight enlargement of the set of test functions necessitated by Proposition \ref{prop:inversion}.

We shortly emphasize that the Hilbert distances of interest in Theorem \ref{thm:sinkhorn} are of course finite, since
\[
d_{\mathcal{G}^{\nu}_{\alpha_2, \tilde{\alpha}_2}}\left(L_{K, \mu}\big(g^{(0)}\big), L_{K, \mu}\left(g_1^*\right)\right) < \infty
\]
by the calculation of the contraction coefficient for $L_{K^\top, \mu}$ from \mbox{Section \ref{sec:integraloperator}}, and immediately obtain
\[
d_{\mathcal{G}^{\mu}_{\alpha, \alpha}}(g^{(n)}, g_1^*) \leq \kappa^{2(n-1)} d_{\mathcal{G}^{\mu}_{\alpha, \alpha}}(g^{(1)}, g_1^*) \leq C \kappa^{2n}.
\] 

There are many different interesting directions that one could explore using Theorem \ref{thm:sinkhorn}, like stability or gradient flows (cf.~\cite{carlier2022lipschitz, DeligiannidisDeBortoliDoucet.21}), obtaining numerically useful bounds, etc. Below, we establish one natural corollary, which transfers the convergence in Hilbert's metric to a more commonly used distance. More concretely, Corollary \ref{cor:Sinkhorn} shows how Theorem \ref{thm:sinkhorn} also yields exponential convergence of the primal iterations of Sinkhorn's algorithm in total variation distance.
To this end, we need to work with a Hilbert metric on the product space $\X \times \X$. We use the metric
\[d_{\X \times \X}((x, y), (\tilde{x}, \tilde{y})) := \max\{d_{\X}(x, \tilde{x}), d_{\X}(y, \tilde{y})\}, ~ x, y, \tilde{x}, \tilde{y} \in \X,\]
and $\bar{B}_r := \{\bar{x} \in \X\times\X : d_{\X\times\X}(\bar{x}, (x_0, x_0)) \leq r\}$.
For Hilbert's metric on the product space, we can for simplicity always consider $\tilde{\alpha} = 0$. For $m > 0$, let
\begin{align*}
\mathcal{F}^{\mu \otimes \nu}_{\alpha, 0} &:= \Big\{f \in L^2(\mu \otimes \nu) : f \geq 0,\\
&\hspace*{3.5cm}\int_{\bar{B}_r^{\com}} f \,d\mu\otimes\nu \leq \alpha(r) \int f \,d\mu\otimes\nu \text{ for all }  r \geq m \Big\},\\
\mathcal{G}^{\mu \otimes \nu}_{\alpha, 0} &:= \left\{g \in L^2(\mu \otimes \nu) : \int fg \,d\mu\otimes\nu \geq 0 \text{ f.a.~}f \in \mathcal{F}^{\mu \otimes\nu}_{\alpha, 0}\right\}.
\end{align*}
The following Lemma shows that we can embed Hilbert metrics from marginal spaces into the product space.
\begin{lemma}\label{lem:embednew}
	Let $\tilde{\alpha}$ be a tail function. For any $f \in \mathcal{F}^{\mu\otimes\nu}_{\alpha, 0}$, we have
	\[
	\int f(\cdot, y) \,\nu(dy) \in \mathcal{F}^{\mu}_{\alpha, \tilde{\alpha}} ~~ \text{ and } ~~ \int f(x, \cdot) \,\mu(dx) \in \mathcal{F}^{\nu}_{\alpha, \tilde{\alpha}}.
	\]
	With the natural embedding of functions mapping from $\mathcal{X}$ into functions mapping from $\X \times \X$ (which are constant in one variable), we have $\mathcal{G}^{\mu}_{\alpha, \tilde{\alpha}} \subseteq \mathcal{G}^{\mu \otimes \nu}_{\alpha, 0}$ and 
	\[
	d_{\mathcal{G}^{\mu\otimes\nu}_{\alpha, 0}}(g, \tilde{g})\leq d_{\mathcal{G}^{\mu}_{\alpha, \tilde{\alpha}}}(g, \tilde{g}) ~~~\text{ for } g, \tilde{g} \in \mathcal{G}^{\mu}_{\alpha, \tilde{\alpha}} ,
	\]
	and also the analogous statement for $\mathcal{G}^{\nu}_{\alpha, \tilde{\alpha}}$ instead of $\mathcal{G}^{\mu}_{\alpha, \tilde{\alpha}}$.
	\begin{proof}
		Let $f \in \mathcal{F}^{\mu\otimes\nu}_{\alpha, 0}$. We only show $f_1 := \int f(\cdot, y) \,\nu(dy) \in \mathcal{F}^{\mu}_{\alpha, \tilde{\alpha}}$. Clearly, $f_1$ is still non-negative, and hence we only have to show the bound for the positive tails. This follows since $B_r^{\com} \times \X \subseteq \bar{B}_r^{\com}$ and hence
		\[
		\int_{B_r^{\com}} f_1 \,d\mu = \int_{B_r^{\com} \times \X} f \,d\mu \otimes \nu \leq \int_{\bar{B}_r^{\com}} f \,d\mu \otimes \nu \leq \alpha(r) \int f \,d\mu \otimes \nu = \alpha(r) \int f_1\,d\mu.
		\]
		The inclusion $\mathcal{G}^{\mu}_{\alpha, \tilde{\alpha}} \subseteq \mathcal{G}^{\mu \otimes \nu}_{\alpha, 0}$ follows from the statement that we have just shown, since for $g \in \mathcal{G}^{\mu}_{\alpha, \tilde{\alpha}}$ and $f \in \mathcal{F}^{\mu\otimes\nu}_{\alpha, 0}$, we have
		\[
		\int gf \,d\mu\otimes\nu = \int g(x) \left(\int f(x, y) \,\nu(dy)\right) \,\mu(dx) \geq 0.
		\]
		Finally, the inequality for the Hilbert metrics follows analogously using Lemma \ref{lem:calcdist}.
	\end{proof}
\end{lemma}

We can now state the exponential convergence result for the primal iterations of Sinkhorn's algorithm.
\begin{corollary}\label{cor:Sinkhorn}
	We make the same assumptions as in Theorem \ref{thm:sinkhorn}.
	Denote by $g_1^*, g_2^*$ the (exponential) optimal dual potentials of the EOT problem with cost $c$, and $g_1^{(n)}, g_2^{(n)}$ the corresponding Sinkhorn iterations. Denote by $\pi^*, \pi^{(n)} \in \mathcal{P}(\X \times \X)$ the primal optimizer and the respective Sinkhorn iterations given by
	\[
	\frac{d\pi^*}{d\mu\otimes\nu}(x, y) = g_1^*(x) g_2^*(y) K(x, y)~~~\text{and}~~~\frac{d\pi^{(n)}}{d\mu\otimes\nu}(x, y) = g_1^{(n)}(x) g_2^{(n)}(y) K(x, y).
	\]
	Then, there exists a constant $C > 0$ and $\rho \in (0, 1)$ such that for all $n \in \mathbb{N}$,
	\[
	\|\pi^{(n)}-\pi^*\|_{TV} \leq C \rho^n.
	\]
	\begin{proof}
		From Theorem \ref{thm:sinkhorn} we know that for $m$ large enough, $g_1^{(n)}$ converges to $g_1^*$ with respect to $d_{\mathcal{G}^{\mu}_{\alpha_1, \alpha_1}}$ and by symmetry also $g_2^{(n)}$ to $g_2^*$ with respect to $d_{\mathcal{G}^\nu_{\alpha_1, \alpha_1}}$, where $\alpha_1(r) = \exp(-r^{p_1})$ with $p_1 = p+\delta/2$. By Lemma \ref{lem:embednew},
		\begin{align*}
		d_{\mathcal{G}^{\mu\otimes\nu}_{\alpha_1, 0}}(g_1^{(n)}, g_1^*) &\leq d_{\mathcal{G}^{\mu}_{\alpha_1, \alpha_1}}(g_1^{(n)}, g_1^*),\\
		d_{\mathcal{G}^{\mu\otimes\nu}_{\alpha_1, 0}}(g_2^{(n)}, g_2^*) &\leq d_{\mathcal{G}^{\nu}_{\alpha_1, \alpha_1}}(g_2^{(n)}, g_2^*).
		\end{align*}
		Next, we define $\alpha_2(r) = \exp(-r^{p+2\delta/3})$. First, we note that by the triangle inequality,
		\begin{align*}
		d_{\mathcal{G}^{\mu\otimes\nu}_{\alpha_2, 0}}(g_1^{(n)} g_2^{(n)}, g_1^* g_2^*) &\leq d_{\mathcal{G}^{\mu\otimes\nu}_{\alpha_2, 0}}(g_1^{(n)} g_2^*, g_1^* g_2^*) + d_{\mathcal{G}^{\mu\otimes\nu}_{\alpha_2, 0}}(g_1^{(n)} g_2^{(n)}, g_1^{(n)} g_2^*).
		\end{align*}
		As in the proof of Theorem \ref{thm:sinkhorn}, with $\kappa(r) := \exp(r^{p+\delta/3})$, for possibly re-scaled versions of $g_1^{(n)}$ and $g_2^*$, we obtain using Lemma \ref{lem:growthnew} that
		\begin{align*}
		1 \leq g_1^{(n)}(x) \leq \kappa(d_{\X}(x_0, x)) \text{ and } 1 \leq g_2^*(x) \leq \kappa(d_{\X}(x_0, x)) ~~~ \text{ for all } x \in B_m^{\com}
		\end{align*}
		for $m$ large enough.
		Thus, using Lemma \ref{lem:multclosedexp}, we have 
		\begin{align*}
		d_{\mathcal{G}^{\mu\otimes\nu}_{\alpha_2, 0}}(g_1^{(n)} g_2^{(n)}, g_1^{(n)} g_2^*) &\leq d_{\mathcal{G}^{\mu\otimes\nu}_{\alpha_1, 0}}(g_2^{(n)}, g_2^*),\\
		d_{\mathcal{G}^{\mu\otimes\nu}_{\alpha_2, 0}}(g_1^{(n)} g_2^*, g_1^* g_2^*) &\leq d_{\mathcal{G}^{\mu\otimes\nu}_{\alpha_1, 0}}(g_1^{(n)}, g_1^*),
		\end{align*}
		and thus
		\begin{align*}
		d_{\mathcal{G}^{\mu\otimes\nu}_{\alpha_2, 0}}(g_1^{(n)} g_2^{(n)}, g_1^* g_2^*) &\leq d_{\mathcal{G}^{\mu}_{\alpha_1, \alpha_1}}(g_1^{(n)}, g_1^*) + d_{\mathcal{G}^{\nu}_{\alpha_1, \alpha_1}}(g_2^{(n)}, g_2^*).
		\end{align*}
		Let $\alpha_3(r) = \exp(-r^{p+3\delta/4})$ for $r \geq 0$. Noting the bounded growth of $K$ through $l$ and \eqref{eq:growthlsink}, another application of Lemma \ref{lem:multclosedexp} yields
		\begin{align}\label{eq:final_ineq}
		\begin{split}
		d_{\mathcal{G}^{\mu\otimes\nu}_{\alpha_3, 0}}(K g_1^{(n)} g_2^{(n)}, K g_1^* g_2^*) &\leq d_{\mathcal{G}^{\mu\otimes\nu}_{\alpha_2, 0}}(g_1^{(n)} g_2^{(n)}, g_1^* g_2^*) \\&\leq d_{\mathcal{G}^{\mu}_{\alpha_1, \alpha_1}}(g_1^{(n)}, g_1^*) + d_{\mathcal{G}^{\nu}_{\alpha_1, \alpha_1}}(g_2^{(n)}, g_2^*).
		\end{split}
		\end{align}
		Thus, using Proposition \ref{prop:relationmetrics} (ii), we find that since
		\[
		\int K g_1^{(n)} g_2^{(n)} \,d\mu\otimes\nu = \int 1 \,d\pi^{(n)} = \int 1 \,d\pi^* = \int K g_1^* g_2^* \,d\mu\otimes\nu,
		\]
		there exists a constant $L > 0$ (arising from Lipschitz continuity of $x \mapsto \exp(x)-1$ when restricted to bounded sets; in our case, say with constant $L/3$) such that
		\[
		\|\pi^{(n)}-\pi^*\|_{TV} = \|Kg_1^{(n)}g_2^{(n)} - K g_1^* g_2^*\|_{L^1(\mu\otimes\nu)} \leq L d_{\mathcal{G}^{\mu\otimes\nu}_{\alpha_3, 0}}(K g_1^{(n)} g_2^{(n)}, K g_1^* g_2^*)
		\]
		which, in view of \eqref{eq:final_ineq}, yields the claim.
	\end{proof}
\end{corollary}

\bmhead{Acknowledgments}

The author thanks Daniel Bartl and Marcel Nutz for helpful comments and discussions. The author is also grateful to the reviewers for their thorough reading and constructive comments.

\section*{Declarations}

The author has no conflicts of interest to declare and is grateful for support by the German Research Foundation through the project 553088969 as well as the Cluster of Excellence ``Machine Learning - New Perspectives for Science'' (EXC 2064/1 number 390727645).

\bibliography{stochfin}

\end{document}